\documentclass[12pt, reqno]{article}

\usepackage{braket}
\usepackage{hyperref}
\usepackage{theoremref}
\usepackage{amsmath}
\usepackage{amsthm}
\usepackage{amssymb}
\usepackage{amsfonts}
\usepackage{color}
\usepackage[all]{xy} \xyoption{all}
 \usepackage{wrapfig}
\usepackage{enumerate}
\usepackage{enumitem}
\usepackage{slashed}
\usepackage{ytableau}

\usepackage{tikz, tikz-3dplot, pgfplots}
\usetikzlibrary{decorations.markings, patterns, decorations.pathmorphing}

\tikzset{->-/.style={decoration={
markings,
mark=at position #1 with {\arrow{>}}},postaction={decorate}}}

\def\CC{\mathbb{C}}

 \def\LL{{\mathbb{L}}}

\def\ZZ{\mathbb{Z}}

\def\TT{{\mathbb{T}}}

\def\aen{\mathfrak{a}}

\def\gen{\mathfrak{g}}

\def\len{\mathfrak{l}}
\def\men{\mathfrak{m}}
\def\nen{\mathfrak{n}}
\def\oen{\mathfrak{o}}
\def\pen{\mathfrak{p}}

\def\sen{\mathfrak{s}}

\def\Cc{\mathcal{C}}

\def\Oc{\mathcal{O}}

\def\Tc{\mathcal{T}}

\def\Vc{\mathcal{V}}

\def\fsl{\mathfrak{sl}}

 \def\act{{\on{act}}}

 \def\be{\begin{equation}}
 \def\bef{\begin{figure}}
 \def\bem{\begin{matrix}}
 
 \def\bpm{\begin{pmatrix}}
  
 \def\btp{\begin{tikzpicture}}

  \def\Cliff{{\on{Cliff}}}
  
 \def\codim{{\on{codim}}}
 
  \def\Coker{\operatorname{Coker}\nolimits}

 \def\del{{\partial}}

  \def\ee{\end{equation}}
 \def\End{\on{{End}}}
 \def\enf{\end{figure}}
 \def\enm{\end{matrix}}
 
 \def\epm{\end{pmatrix}}
 \def\eps{{\varepsilon}}
 \def\etp{\end{tikzpicture}}
 \def\even{{\on{even}}}

 \def\fd{{\on{fd}}}
 \def\FMod{{\on{FMod}}}

 \def\gl{{\gen\len}}
 \def\GL{{\on{GL}}}

  \def\Hom{\operatorname{Hom}\nolimits}

   \def\hra{\hookrightarrow}

 \def\Id{\operatorname{Id}\nolimits}
 \def\Ind{{\on{Ind}}}

 \def\kappa{\varkappa}
   \def\Ker{\operatorname{Ker}\nolimits}

  \def\Lie{{\on{Lie}}}

   \def\lra{\longrightarrow}

  \def\Mod{{\text{-}\on{Mod}}}

\def\o{{\mathfrak{o}}}

\def\odd{{\on{odd}}}
 
\def\on{\operatorname}
\def\ol{\overline}

\def\oo{{\infty}}

\def\osp{{\mathfrak{osp}}}

 \def\PF{{\on{PF}}}
 \def\phi{{\varphi}}

\def\Rep{{\on{Rep}}}

  \def\sgn{{\on{sgn}}}

  \def\ssl{\mathfrak{{sl}}}

  \def\sp{{\mathfrak{sp}}} 
  \def\SPF{{\on{SPF}}}
  \def\Spec{{\on{Spec}}}
  
  \def\ssl{{\sen\len}}
  
  \def\str{{\on{str}}}

 \def\SVect{{\on{SVect}}}

 \def\Tor{{\on{Tor}}}
 \def\tr{{\on{tr}}}

\def\ul{\underline}

\def\Vect{\on{Vect}}

\def\wh{ \widehat}
\def\wt{\widetilde}

\def\YD{{\on{YD}}}

\def\Z2YD{{{1\over 2}\ZZ\hskip -0.07cm\on{YD}}}

\def\0{{\ol{0}}}
\def\1{{\ol{1}}}
\def\b1{{{\mathbf{1}}}}
\def\(({(\hskip -1mm (}
\def\)){)\hskip -1mm )}
\def\-{{\,\, \setminus} \,\,}
\def \= {{\, \simeq \,}}
 \def\be{\begin{equation}}
\def\ee{\end{equation}}
\def\ed{\end{document}}

\newtheorem{thm}[equation]{Theorem}

\newtheorem{lem}[equation]{Lemma}
\newtheorem{prop}[equation]{Proposition}

\theoremstyle{definition}

\theoremstyle{remark}

\newtheorem{rem}[equation]{Remark}

\newtheorem{ex}[equation]{Example}
\newtheorem{exas}[equation]{Examples}

\numberwithin{itemcounter}{section}
\numberwithin{equation}{section}
\title{ An explicit BGG-type resolution of the spin/oscillator module for
$\osp(2m+1|2n)$ }

\author{Mikhail Kapranov, Vera Serganova}

\DeclareUnicodeCharacter{2297}{\ifmmode\otimes\else{$\otimes$}\fi}
\DeclareUnicodeCharacter{2194}{\ifmmode\leftrightarrow\else{$\leftrightarrow$}\fi}

\begin{document}
\maketitle
 
 \abstract{ We construct an explicit resolution of the spin/oscillator representation
 of the Lie superalgebra $\osp(2m+1|2n)$ by parabolic Verma modules. 
 The resolution appears as the
 Chevalley-Eilenberg chain complex of the Lie algebra of quadratic operators
 in the universal Heisenberg algebra of the vector superspace $\CC^{m|n}$. 
 In particular, its differentials are given by unified linear algebra
   formulas and all have degree $2$
 with respect to  the standard grading of Verma modules. 
  }
  
   \tableofcontents

  \section*{Introduction}
  
  \paragraph{Motivation.} 
  The classic Bernstein-Gelfand-Gelfand (BGG) resolution,
   see, e.g.,  \cite{humphreys},  is a fundamental
  feature of representation theory of a complex semisimple Lie algebra $\gen$. It 
  resolves a finite-dimensional irreducible $\gen$-module $M$ via a complex
  of induced Verma modules labelled by elements of the Weyl group $W$
  (or, for  parabolic BGG resolutions,  of the set of parabolic representatives in  $W$). The identical shape of the resolutions for different  such $M$ is
  explained by the structure of the category $\Oc^\gen$ of $\gen$-modules
  generated by highest weight vectors which splits into blocks labelled
  by $W$-orbits on weights. The blocks containing different  finite-dimensional irreducibles
  are all equivalent and commonly referred to as the {\em regular blocks}. 
  
  \vskip .2cm
  
  There are interesting representations in $\Oc$ which do not belong to regular blocks, the most
  immediate being the {\em Weil } (or {\em oscillator}) representation of $\gen=\sp(2n)$ in the polynomial space
  $\CC[x]=\CC[x_1,\cdots, x_n]$. Its irreducible components are the spaces $\CC[x]_\0$ and $\CC[x]_\1$ of
  even and odd polynomials. Their parabolic BGG-type resolutions can be obtained via the Enright-Shelton equivalence
  \cite{enright-shelton} 
   identifying the corresponding (irregular) block of the parabolic $\Oc^{\sp(2n)}$ with the sum of two regular blocks of   
  a parabolic $\Oc^{\o(2n+1)}$.   This does not give a good handle of the differentials in the resolution. For instance, 
restriction to the enveloping of the (abelian) nilpotent radical of the maximal parabolic in $\sp(2n)$
gives the J\'oseffiak-Pragacz-Weyman resolution   of the ideal of matrices of rank $\leq 1$  in  the coordinate ring of
  the space of symmetric matrices  \cite[Th.3.19]{weyman-resol},
  see also \cite[Th.6.3.1]{weyman-book}. In this case, as explained in {\em loc. cit.},  the matrix elements
   of the differentials are hard to find explicitly and  their
  degrees can be arbitrary large.

  \paragraph{Our result.} It turns out that passing to Lie superalgebras simplifies things. In this paper we consider the
  Lie superalgebra $\gen = \osp(2m+1|2n)$ and its spin/oscillator representation in the space $\Lambda^\bullet(L)$,
  where $L=\CC^{m|n}$ is a maximal isotropic space in the orthogonal  superspace $\CC^{2m+1|2n}$. 
  For $n=0$ we get the spin representation of $\o(2m+1)$ and for $n=0$ we get $\CC[x]$ as a module
  over $\osp(1|2n)$,  uniting in one irreducible  the two summands of the oscillator representation of $\sp(2n)$, cf.  \cite{malcom, kapranov}. 
  
   Our main result, Theorem \ref{thm:existence},  gives a resolution
   \be\label{eq:res-intro}
 \cdots \lra \Vc^{-2} \lra \Vc^{-1} \lra \Vc^0
  \to \Lambda^\bullet(L) \to 0,
\ee
where $\Vc^{-p} = V_{\Lambda^p(S^2(L))\otimes \CC_{-1/2}}$ is the parabolic Verma module corresponding to the (reducible, in general)   representation $\Lambda^p(S^2(L))\otimes \CC_{-1/2}$
of the Levi  subalgebra $\men =   \gl(L) = \gl(m|n)$. Here $\CC_{-1/2}$ is the character $x\mapsto
-{1\over 2} \str(x)$ of $\men$. A decomposition of $\Vc^{-p}$
into more standard-looking Verma modules can be obtained from the well known
 irreducible decomposition (cf. Proposition \ref {prop:lambda-s2} below)
\[
\Lambda^p(S^2(L)) \= \bigoplus_{\alpha\in Q_{-1}(2p)} \Sigma^\alpha(L).
\]
Here $\Sigma^\alpha$ is the Schur functor and  $Q_{-1}(2p)$ is the set of Young diagrams $\alpha$ having 
 $2p$ boxes with  the property that 
  each arm of $\alpha$ 
is   by one  longer than the corresponding leg. 

\vskip .2cm

For $n=0$ this gives the standard BGG resolution of the spinor module for $\o(2m+1)$
and the identification of the terms with $V_{\Lambda^p(S^2(L))\otimes \CC_{-1/2}}$
was noted in \cite{SS-spin}. For $m=0$ the existence of such resolution 
was proved
in \cite{kapranov} using the Enright-Shelton equivalence mentioned above together
with the Gorelik equivalence \cite {gorelik-typical} relating the categories of representations
of $\gen$ and its even part $\sp(2n)$.  In general, the resolution is infinite
and we do not know a parametrization of its summands by elements of some Weyl-type group. 
The statements implied by our resolution
at the level of characters (and even of the Kostant homology) were known before
(\cite {chen-kwon-wang}).

 \paragraph{Our method: explicit construction as the Chevalley-Eilenberg complex.}  Unlike the usual approaches based on combinatorics of Weyl groups or topology of
Schubert varieties, our construction of BGG-type resolutions proceeds by pure linear algebra.

  We first work over the (opposite) nilpotent radical $\ol\nen = L\oplus \Lambda^2(L)$
  so that Verma modules are free over $U=U(\ol\nen)$. This $U$ is the universal
  Heisenberg algebra over its center $Z= S^\bullet(\Lambda^2L)$ and therefore
  the $Z$-module $\aen = S^2(L)\cdot Z$ of {\em quadratic Heisenberg operators} is a Lie algebra over $Z$ via the classic Moyal formula. Our resolution (considered as
  a complex of $Z$-modules) is nothing but the Chevalley-Eilenberg chain complex
  of $\aen$ with coefficients in $U$. In particular, its differentials are given by
  the universal formulas and all have degree $2$ (with respect to the natural grading
  of Verma modules). 
  
  We then prove Theorem \ref{thm:exactness} which says 
  that this indeed gives a resolution of $\Lambda^\bullet(L)$, so far by $U$-modules.
   Next, we check
  (Theorem   \ref{thm:linearity}) that after
  tensoring with $\CC_{-1/2}$, the Chevalley-Eilenberg complex becomes a
  complex of $\gen$-modules.

  Thus even for the spin representation of $\o(2m+1)$ when the standard BGG resolution is
  available, our approach gives something new: a unified linear algebra construction
  of the differentials. 
  
  \paragraph{The tensor category context.}  Our constructions are functorial in the vector space $L=\CC^{m|n}$.
  The algebras $\ol\nen, \aen, U, Z$, the module $\Lambda^\bullet(L)$ and the complex $\Cc^\bullet$
   are obtained by applying to $L$ some polynomial functors
  (direct sums of Schur functors) and the differentials in $\Cc^\bullet$  have matrix elements which are
   natural morphisms of Schur functors.
  The Lie algebra $\gen$ is obtained similarly but by using the dual space $L^*$ as well. 
  This suggests that   our constructions can be performed   in more general  tensor categories. 
  This  is indeed so and we explain this in
   \S \ref {sec:tensor}. The right context is that of a $\CC$-linear symmetric monoidal abelian category $\Tc$ with
   a {\em bilinear pair}, i.e., a pair of objects $L, L_*$ and a pairing $\omega: L_*\otimes L\to \b1$ with values in the unit object.
 The universal example $\TT_\ssl$ of such category, studied in    \cite{DPS, SS-stable}, consists of certain representations of 
 $\ssl(\oo) = \varinjlim \, \ssl(n)$
 or of similarly defined $\GL(\oo)$.   We use this universal example to extend our results  to the more general categorical
 context just described.

  \paragraph{Acknowledgements.} M. K. would like to thank T. Melia for bringing the
  paper \cite{SS-spin} to his attention. His work was supported by the
   World Premier International Research Center Initiative (WPI), MEXT, Japan.

    \vfill\eject

 \section {  Generalities on Young digrams and super-Schur functors.}
 
 \paragraph{Conventions.} \label{par:conventions}
 We work over the field $\CC$ of complex numbers. 
 We denote by  $\Vect$ and $\SVect$ the category of vector spaces resp. vector superspaces over $\CC$.
An element of $\SVect$ is written as $V=V_\0 \oplus V_\1$. The parity of a homogeneous  element $v\in V$
is denoted by $|v|$. Let $\Vect^\fd$ resp. $\SVect^\fd$ be the category of finite-dimensional vector spaces 
resp. superspaces and $\Vect^\fd, \SVect^\fd$ the full subcategories of finite-dimensional spaces there. 

\vskip .2cm

We consider $\SVect$ and $\SVect^\fd$ as $\CC$-linear {\em supercategories},
i.e., as categories enriched in $\SVect$, not just $\Vect$. That is, for any two vector
superspaces $V, V'$ we consider the superspace $\Hom_\CC(V, V')$  
of all $\CC$-linear maps $V\to V'$ (with its $\ZZ/2$-grading induced from that on 
$V$ and $V'$) 
while morphisms from $V$ to $V'$ in the usual (non-enriched) sense are elements
of the even part $\Hom_\CC^\0(V, V')$. 

\vskip .2cm

 For any associative superalgebra $A$ let $A_\Lie$ be $A$
 considered as a Lie superalgebra with bracket (super-commutator)
   \[
  [a,b]= ab- (-1)^{|a|\cdot |b|} ba.
  \]
We will also use the notation
\be\label{eq:super-anticom}
\{a,b\} = ab +  (-1)^{|a|\cdot |b|} ba, \quad a,b\in A
\ee
for the super-symmetrized product (``super-anticommutator'') in $A$.

 \paragraph{Representations of $\gl(m|n)$.}  Let $\gl(m|n) = \End(\CC^{m|n})_\Lie$ be the general linear superalgebra. 
 Dominant  (not necessarily integral) weights of
$ \gl(m|n)$ are sequences $(\lambda|\mu) = (\lambda_1, \cdots, \lambda_m|\mu_1,\cdots,\mu_n)$
of complex numbers such that each $\lambda_i-\lambda_{i+1}$ and each
$\mu_j-\mu_{j+1}$ is an element of $\ZZ_{\geq 0}$. 
Such sequences label 
  irreducible finite-dimensional representations of $\gl(m|n)$,
  see \cite[Ch. 10]{musson}. 
  We denote by $L_{(\lambda|\mu)}$ the representation corresponding to
  $(\lambda|\mu)$. For each $\eps\in \CC$ we will write
  \be\label{eq:C_eps}
  \CC_\eps\,  := \, L_{(\eps, \cdots, \eps | - \eps, \cdots, -\eps)}. 
  \ee
 This is a $1$-dimensional representation corresponding to the character
 $x\mapsto \eps \on{str}(x)$ of $\gl(m|n)$.

 \paragraph {Young diagrams.} 
We understand 
 {\em Young diagrams}  as {\em partitions}, i.e., as sequences
 $\alpha = (\alpha_1\geq\cdots\geq \alpha_r\geq 0)$ of indeterminate length
 consisting of integers and which are considered up to adding a string of $0$'s
 at the end. Explicitly, the nonzero $\alpha_i$ are considered as the numbers of cells
(boxes)  in the  rows of the Young diagram. We write $|\alpha|=\sum \alpha_i$ for the
 number of boxes in $\alpha$. By $\alpha^\vee$ we denote the transpose Young diagram
 to $\alpha$,
 whose rows are  the columns of $\alpha$. Let $\YD$ be the set of all Young diagrams.
 
 \vskip .2cm
 
  We will also use the {\em Frobenius notation} for Young diagrams
  \cite[\S 1.1.2]{weyman-book}, writing
 $\alpha = (a_1,\cdots, a_s |b_1,\cdots, b_s)$ 
 with $a_1>\cdots > a_s > 0$ and $b_1>\cdots > b_s > 0$. 
 Here $a_i$ is the number of boxes in the $i$th {\em arm} of $\alpha$,
 i.e., in the segment  to the right of the $i$th diagonal box (including it) 
 or, in other words, $a_i = \alpha_i-i+1$ as long as this number is $\geq 0$. 
  Similarly, $b_i$ is the number of boxes in the $i$th {\em leg}
 of $\alpha$, i.e., in the segment  below the $i$th
 diagonal box (including it)  or, in other words, $b_i = \alpha^*_i-i+1$ as long
 as  this number is $\geq 0$.  Thus,
 \[
 (a_1,\cdots, a_s |b_1,\cdots, b_s)^\vee = (b_1,\cdots, b_s | a_1,\cdots, a_s).
 \]
 
 \paragraph {Polynomial representations. (Super-) Schur functors.} 
 We will be specially interested in representations of $\gl(m|n)$ which integrate
 to {\em polynomial representations} of the algebraic group $\GL(m|n)$. 
 It is known that the category of polynomial representations of $\GL(m|n)$ is semisimple and
 its irreducible objects are labelled by certain Young digrams
 see \cite {berele-regev} \cite[\S 3.2]{cheng-wang-book}.

\begin{figure}[h]
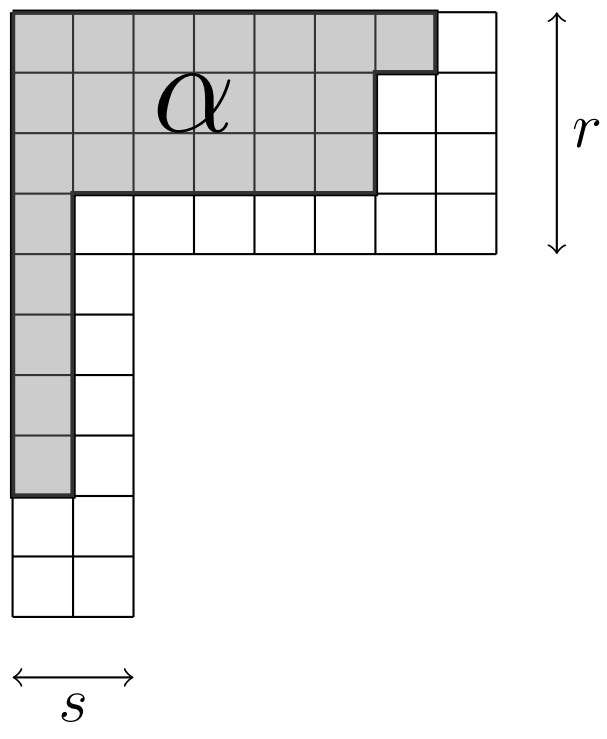
 
   \centering
   \btp[scale=0.4]
   \draw (-2,0) -- (6,0); 
      \draw (-2,1) -- (6,1); 
     \draw (-2,2) -- (6,2);  
      \draw (-2,3) -- (6,3); 
  \draw (-2,4) -- (6,4);

  \draw (-2, -1) -- (0, -1);    
    \draw (-2, -2) -- (0, -2);    
      \draw (-2, -3) -- (0, -3);    
        \draw (-2, -4) -- (0, -4);    
          \draw (-2, -5) -- (0, -5);    
            \draw (-2, -6) -- (0, -6);    
            
\draw (-2, 4) -- (-2, -6); 
\draw (-1, 4) -- (-1, -6); 
\draw (0, 4) -- (0, -6); 

\draw (1,4) -- (1,0); 
\draw (2,4) -- (2,0); 
\draw (3,4) -- (3,0); 
\draw (4,4) -- (4,0); 
\draw (5,4) -- (5,0); 
\draw (6,4) -- (6,0); 
   
 \draw [<->] (7,0) -- (7,4);     
  \draw [<->]  (-2, -7) -- (0, -7); 
  
  \node at (7.5,2) {$r$};
  \node at (-1, -7.5) {$s$};  
  
   \draw[ line width=1] (-2,4) -- (5,4) -- (5,3) -- (4,3) -- (4,1) -- (-1,1) --(-1, -4) --(-2, -4) -- (-2, 4); 
  
  \draw[fill, color={gray}, opacity = 0.4 ] (-2,4) -- (5,4) -- (5,3) -- (4,3) -- (4,1) -- (-1,1) --(-1, -4) --(-2, -4) -- (-2, 4); 
  
 \node at (1, 2.5){\huge$\alpha$}; 
   \etp
   \caption{The hook condition for $\alpha\in \YD(r,s)$.}
   \label{fig:hook-cond}
   \end{figure}

 More precisely, 
 a  Young diagram $\alpha$ gives rise to  the {\em (super-) Schur functor} $\Sigma^\alpha: \SVect \to\SVect$.
  As in the usual (non-super) case
 \cite{weyman-book},
for a vector superspace $V$ the space $\Sigma^\alpha (V)$ can be defined
 inside the tensor power $V^{\otimes |\alpha|}$ as the image of the Young symmetrizer
 associated to $\alpha$, see \cite {berele-regev}. Note the standard identifications
 \[
 \Sigma^k(V) = S^k(V), \,\,\,\Sigma^{1,\cdots,  1} (V) = \Lambda^k(V)
 \text{ (with $k$ occurrences of $1$)}, \,\,\, \Sigma^\alpha(\Pi V) = \Sigma^{\alpha^\vee}(V).
 \]
 For any integers $r,s\geq 0$ let
 \[
 \YD(r|s) = \bigl\{ \alpha \in \YD \,\bigl| \alpha_{r+1} \leq s\bigr\}. 
 \]
 That is, $\alpha$ is contained in the union of its first  $r$ rows and  its
 first $s$ columns, see Fig. \ref {fig:hook-cond}. 
 The condition $\alpha_{r+1}\leq s$ is known as the {\em hook condition}. 
 It is known \cite{berele-regev, cheng-wang-book}
   that for $\dim(V)=(r|s)$ the space $\Sigma^\alpha (V)$
  is nonzero, if
 and only if $\alpha\in \YD(r|s)$. 

 \vskip .2cm
 
 Putting $V=L=\CC^{m|n}$ we have \cite{berele-regev, cheng-wang-book}  that the  $\Sigma^\alpha(L)$, $\alpha\in \YD(m|n)$,
 are precisely all the
 irreducible polynomial representations of $\GL(m|n)$.
  We will denote these representations
 simply by $\Sigma^\alpha$. 
 The description of
 the  highest weight
 $\alpha^\natural$ of $\Sigma^\alpha$ as a representation of
 $\gl(m|n)$ can be found in \cite {cheng-wang-book}, see 
 Eq.(2.4) resp.  Thm. 3.11  there for the definition resp.  the statement and proof. 
 We will not need its precise form. 
 
  \vskip .2cm
  
 Let $l\in \ZZ$ and $p\in \ZZ_{\geq 0}$. 
  As in \cite[\S2.3]{weyman-book}, we denote by $Q_{l}(2p)$ the set
  of Young diagrams $\alpha$ with $|\alpha|=2p$ such that in the
  Frobenius notation $\alpha = (a_1,\cdots, a_s| b_1,\cdots, b_s)$
  we have $b_i = a_i + l$. 
  
  \begin {prop}\label{prop:lambda-s2}
  We have a natural  isomorphism (i.e., an isomorphism of polynomial 
  representations of $\GL(m|n)$)
  \[
  \Lambda^p(S^2(L)) \= \bigoplus_{\alpha\in Q_{-1}(2p)} \Sigma^\alpha (L). 
  \]
  \end{prop}
  \noindent By the above, nonzero summands in this decomposition correspond to
  $\alpha\in Q_{-1}(2p) \cap \YD(m|n)$. 
  
  \vskip .2cm
  
  \noindent {\sl Proof:} For $L = \CC^m$ purely even this is 
  \cite[Prop.2.3.9]{weyman-book}. In general, as Schur functors
  are defined in terms of representations of symmetric groups, the statement
  for $L=\CC^m$ for all $m$  implies an identification of representations of $S_{2p}$
  corresponding to the LHS and the RHS of the proposed general identification.
 More precisely, by Sergeev--Schur--Weyl duality \cite{sergeev}
      we have that
      \[
      [\Lambda^p(S^2L):\Sigma^\alpha(L)]=[\Ind^{S_{2p}}_{B_p }\sgn:\Gamma_\alpha],
      \]
    where $B_p=S_p\ltimes(\mathbb Z/2)^p$ (the Weyl group for $B_p$ root system),  $\sgn$ is inflated from the sign representation of $S_p$, and $\Gamma_\alpha$ the irreducible representation of $S_{2p}$ associated with $\alpha$.
      Therefore it implies this general identification.
  \qed
  
  \paragraph{Polynomial functors on $\SVect^\fd$.}

  Extending the standard (purely even) concept, see \cite[App. A]{macdonald}, by a 
   {\em (super) polynomial functor} we mean an (super-enriched, in the sense 
   of \S \ref {par:conventions}) functor
  $F: \SVect^\fd \to \SVect$ such that for any $V, V'\in \SVect^\fd$
  the map from the superspace $\Hom_\CC(V, V')$ to the
  superspace $\Hom_\CC(F(V), F(V'))$ induced by $F$, is a (super) polynomial map. 
  Let $\SPF$ be the category of polynomial functors on $\SVect^\fd$. It has a tensor
  structure given by objectswise tensor product of functors: $(F'\otimes F'')(V) =
  F(V)\otimes F''(V)$. 
  
  \vskip .2cm

 Any $F\in\SPF$ induces a (possibly infinite-dimensional) polynomial representation
  of  $\GL(V)\=  \GL(r|s)$ for any $V\= \CC^{r|s}\in\SVect^\fd$. 
  Because of semi-simplicity of the category of polynomial representations of $\GL(r|s)$,
  this implies, as in \cite[App.A]{macdonald},
   that $F$ is a (possibly infinite) direct sum of Schur functors:
  $F\= \bigoplus_{\alpha\in YD} C_\alpha \otimes \Sigma^\alpha$ for some
  (possibly infinite-dimensional) vector superspaces $C_\alpha$. 
  We say that $F$ is {\em even}, if all the $C_\alpha$ are purely even. 
  Let $\SPF^\0\subset \SPF$ be the full subcategory of even super-polynomial functors. 
  An alternative description of polynomial functors is as follows,
  cf. \cite{getzler, joyal, macdonald}. 
  
  \vskip .2cm

Let  $S_d$ be  the symmetric group on $d$ letters. By an {\em $S$-(super)module}
we will mean a collection $M= (M_d)_{d\geq 0}$ of vector superspaces
with $S_d$ acting on $M_d$ for each $d$. Let $S\Mod$ denote the category
of $S$-modules.  It has a tensor structure denote $*$.  That is,
 for $M', M''\in S\Mod$ the
$S$-module $M'*M''$ is defined by
\be
  (M' *M'')_d = \bigoplus_{d'+d''=d} \Ind_{S_{d'}\times S_{d''}}^{S_d} (M'_{d'}\otimes M''_{d''}).
\ee
The following is then a reformulation of the super  version 
\cite{berele-regev, cheng-wang-book}  of the Schur-Weyl duality  between representations
of symmetric groups and polynomial representations of the $\GL(r)$. 

\begin{prop}\label{prop:pol-fun}
  Associating to an  $S$-module $M$  the  polynomial functor $F_M$ by
  \[
  F_M(V) = \bigoplus_d \, (V^{\otimes d}\otimes M_d)_{S_d}
  \]
  establishes an equivalence   of tensor categories $(S\Mod, *) \to (\SPF, \otimes)$. 
   \qed
 \end{prop}

  Let $S\Mod^\0\subset S\Mod$ be the full subcategory of even $S$-modules
  i.e., of $M$ such that each $M_d$ is purely even. Let also $\PF$ be the
  category of polynomial functors $\Vect^\fd\to\Vect$ in the classical (even)
  sense  \cite{joyal}. 
  Then the classical (even) version of Proposition \ref {prop:pol-fun} given an 
  equivance 
  $
  S\Mod^\0 \= \PF 
  $. On the other hand,  the equivalence of Proposition \ref {prop:pol-fun}
   identifies  $S\Mod^\0$  with
  $\SPF^\0$, the category of even polynomial functors on $\SVect^\fd$.
  Therefore $\SPF^\0 \= \PF$. 
 Now, each $S$-module $M$ can be written as  direct sum of $M' \oplus \Pi M''$
  where $M', M''$ are even and so $\SPF$ is identified with the direct sum of
  categories $\PF = \SPF^\0$ and $\Pi(\PF)$. 
   This implies the following.

  \begin{prop}\label{prop:exactness-even}
  Let $F\to G\to H$ be a sequence of  polynomial functors   in $\SPF$  and natural transformations.
  Then the following are equivalent:
  \begin{itemize}
  \item[(i)] The sequence $F\to G\to H$ is exact, i.e., for any $V\in\SVect^\fd$
  the sequence $F(V)\to G(V)\to H(V)$ is exact.
  
  \item[(ii)]For any purely even vector space $V\in \Vect^\fd$ the sequence $F(V)\to G(V)\to H(V)$
  is exact.
  
  \item[(iii)] For any purely odd vector space $V= \Pi(V'), V'\in\Vect^\fd$ 
  the sequence $F(V)\to G(V)\to H(V)$
  is exact.

  \end{itemize}
  \qed
  \end{prop}
  
  \begin{rem}\label{rem:vanishing-even}
  Similarly, we can verify vanishing 
  of a morphism $\phi: F \to G$ of polynomial functors  in $\SPF$ on purely even or
  purely odd spaces only. In fact, this is a consequence of Proposition 
  \ref {prop:exactness-even}, as vanishing of $\phi$
   is equivalent to exactness of the sequence $F\buildrel = \over\to F\buildrel\phi\over \to G \buildrel=\over\to  G$. 
  \end{rem}

  \vfill\eject

  \section{Construction of the resolution: statement of results}\label{const-resul}
  
  \paragraph{Background on $\osp(2m+1|2n)$.}
  
  Let $V$ be a finite-dimensional vector superspace with
  a nondegenerate scalar product (a symmetric, in the super sense,    bilinear
  form)  $(-,-)$. Denote by $\o(V)$ the orthogonal Lie superalgebra
  of $(V, (-,-))$, consisting of linear operators $A: V\to V$ such that
  $A^t = -A$. Thus, if $\dim(V)=(r|s)$, then $\o(V) = \osp(r|s)$
  is the orthosymplectic Lie superalgebra in the standard sense. 
  
  \vskip .2cm
  
  Identifying the source of $A\in \o(V)$ with $V^*$ via $(-,-)$,
  we can view $A$ as an element of $\Lambda^2(V)$. In other
  words, we have an identification $ \o(V)\= \Lambda^2(V)$ which is
  in fact an identification of $\o(V)$-modules, where the action
  on $\o(V)$ is the adjoint one. More explicitly, the wedge product
  $v_1\wedge v_2$ acts on $V$ by
  \[
  (v_1\wedge v_2) (v) = (v_2, v) \cdot v_1 - (-1)^{|v_1|\cdot |v_2|}
  (v_1, v) \cdot v_2. 
  \]
  
  \vskip .2cm
  
    Let us now start from  a vector superspace $L=\CC^{m|n}$
   of dimension $m|n$. The space 
  $E=L\oplus L^*$ has the hyperbolic scalar product 
  \[
  \bigl((v,f), (v',f')\bigr) =  f(v') + (-1)^{|f'|\cdot  |v| } f'(v). 
  \]
 Consider further the superspace $\ol E = E\oplus \CC$
 of dimension $(2m+1|2n)$ with the direct
 sum of the above scalar product on $E$ and the
 scalar product $xx'$ on $\CC$,
 i.e.,
 \[
 \bigl( (v,f,x), (v', f', x')\bigr) =  \bigl((v,f), (v',f')\bigr) + xx'. 
 \]
 Let 
\[
 \gen = \o(\ol E) = \osp(2m+1|2n)
 \]
 be the orthogonal Lie superalgebra of this superspace.
 By the above, we can identify
 \[
 \gen = \Lambda^2(E\oplus \CC) = E\oplus \Lambda^2(E), \quad \Lambda^2(E)=\o(E). 
 \]
 This decomposition can be seen as an additional 
 (with respect to the standard super-grading) $\ZZ/2$-grading
 of the Lie superalgebra $\gen$, with $E$ sitting in degree $\1$
 and $\Lambda^2(E)$ in degree $\0$. 
 The degree $\0$ part, i.e., $\Lambda^2(E)$,
 is a sub-superalgebra, namely $\o(E)$;  the 
 part of the bracket  $E\otimes E\to\Lambda^2(E)$
 is the wedge product;  the action of $\Lambda^2(E)=\o(E)$
 on $E$  given by the commutator  is the standard action of $\o(E)$
 on $E$. 
 
 \vskip .2cm
 
 We further  have the parabolic decomposition
 \be\label{eq:g-notation}
 \begin{gathered}
 \gen = \Lambda^2(\ol E) = \Lambda^2 (L\oplus L^* \oplus \CC) = \ol\nen \oplus \men \oplus
 \nen,
 \\
 \ol\nen := L\oplus \Lambda^2(L), \quad \men := L\otimes L^* \= \gl(m|n),
 \quad \nen := L^* \oplus \Lambda^2(L^*),
 \\
 \text{ so that}
 \\
  \ol \pen := \men + \ol\nen, \quad \pen := \men \oplus\nen
 \end{gathered}
 \ee
 are 
 parabolic subalgebras with the same Levi $\men$ and radicals
 $\ol\nen$ and $\nen$. 
 The subalgebra $\ol\nen$ is the stabilizer of the isotropic subspace
 $L\subset \ol E$ and $\nen$ is the stabilizer of $L^*\subset \ol E$.
 
 \vskip .2cm
 
 As Lie superalgebras, $\ol\nen$ and $\nen$ are free nilpotent
 of degree $2$. That is, the only nonzero component of the bracket
 in $\ol \nen$ is the map $L\otimes L\to\Lambda^2(L)$
 which is the wedge product, and similarly for $\nen$. 
 
 \paragraph{The spin/oscillator representation of $\osp(2m+1|2n)$.}

 Consider the   Clifford superalgebra $\Cliff(E)$, $E=L\oplus L^*$ which
 can be seen as the algebra of (super) differential operators in
 the  exterior algebra
 $\Lambda^\bullet(L) = S^\bullet(\Pi L)$. 
  Explicitly, let $y_1,\cdots, y_m$ be a basis
 of $L_\0$ and $y_{m+1},\cdots, y_{m+n}$ be a basis of $L_\1$. 
 Let $z_1,\cdots, z_{m+n}$ be the dual basis of $L^*$. 
 We write $\Lambda^\bullet(L)$ as 
 $\Lambda^\bullet[y_1,\cdots, y_{m+n}]$. Then $y_i\in\Cliff(E)$ acts in
 $\Lambda^\bullet(L)$ as left multiplication by $y_i$, while
 $z_i$ acts by the (super) derivative $\del/\del y_i$. 
 
  \vskip .2cm
  
 The Lie superalgebra
  $\o(E)$ is realized inside $\Cliff(E)_\Lie$ as the subspace
  $\Lambda^2(E)$ of {\em quadratic Clifford operators}, i.e.,
  of linear combinations of  the commutators $[a,b]$ with $a,b\in E$.
  That is, $\o(E)=[E,E]\subset \Cliff(E)_\Lie$. This makes
  $\Lambda^\bullet(L)$ into an $\o(E)$-module. As such,
  this module is reducible, splitting as
  \[
  \Lambda^\bullet(L) \, = \, \Lambda^\even (L) \oplus \Lambda^\odd (L),
  \quad \Lambda^\even(L) := \bigoplus_k \Lambda^{2k}(L), \,\,\,
  \Lambda^\odd(L) := \bigoplus_k \Lambda^{2k+1}(L). 
  \]
  Note that the summands $\Lambda^\even(L)$ and $\Lambda^\odd(L)$
  are, in general, not the same as $\Lambda^\bullet(L)_\0$ and
  $\Lambda^\bullet(L)_\1$,  the even and odd parts of
  $\Lambda^\bullet(L)$ as a vector super space.  These  summands,
  which are vector superspaces in their own right, 
  are the two {\em spin/oscillator representations} of 
  $\o(E) =\osp(2m|2n)$. 
  
  \vskip .2cm
  
  Next, we extend the action of $\o(E)=\Lambda^2(E)$ on $\Lambda^\bullet(L)$
  to an action of  the superalgebra
  $\o(\ol E) =  E\oplus \Lambda^2(E)$. For this, we define
  the action of $E=L\oplus L^*$ to be the restriction of the action of $\Cliff(E)$
  as defined above. The following is then well known, cf. \cite {hartwig-S}.
  
  \begin{prop}\label{prop:spin/weil}
 (a)  This rule   defines an action of $\gen = \o(\ol E) = \osp(2m+1|2n)$ on 
 $\Lambda^\bullet(L)$.
 
 \vskip .2cm
 
(b)  The subspace $\CC\cdot 1\subset \Lambda^\bullet(L)$ is invariant under
 $\men=\gl(m|n)$ which acts there via the character $\CC_{-1/2}$. 
  \qed
  \end{prop}

The space $\Lambda^\bullet(L)$ with this action is known as the
{\em spin/oscillator representation}  of $\gen = \osp(2m+1|2n)$. 
It is irreducible since the action of $E$ generates $\Cliff(E)$ whose
action on $\Lambda^\bullet(L)$ is irreducible. 

\begin{exas}\label{exas:spiun/osc}
(a) Let $n=0$. Then $\gen=\o(2m+1)$ is the usual orthogonal Lie  (non-super) algebra
of degree $2m+1$. In this case $\Lambda^\bullet(L)$, $L=\CC^m$,  is the unique spinor
representation of $\o(2m+1)$. After restriction to $\o(2m)$ it splits into the sum
of the two spinor representation of $\o(2m)$. 

\vskip .2cm 

(b) Let $m=0$. Then  $\gen = \osp(1|2n)$. In this case $\Lambda^\bullet(L) = S^\bullet
(\CC^n)$ is the oscillator (Weil) representation of $\osp(1|2n)$, see 
\cite{malcom, kapranov}.  After restriction to $\o(E)=\sp(2n)$ it splits into the
sum of the classical even and odd Weil representations $S^\even(\CC^n)$ and
$S^\odd(\CC^n)$ which form (in this case only) the 
$\ZZ/2$-grading giving  the structure of  a vector superspace on $\Lambda^\bullet(L)$,
cf. \cite{kapranov}.

\end{exas}

\vfill\eject\vfill\eject

\paragraph{The parabolic category $\Oc$ and the BGG resolution.} 
Let $\Oc=\Oc^{\gen, \pen}$ be the category of $\gen$-modules 
  $V$ which are:
\begin{itemize}
\item Finitely generated over $U(\gen)$.

\item Considered as a $U(\men)$-module, split into a direct sum of finite-dimensional
simple modules  with semisimple action of the center of $\men$.

\item Locally $U(\nen$)-finite, i.e., the $\nen$-module $U(\nen)v$ is finite-dimensional
for any $v\in V$. 
\end{itemize} 

\noindent For example, $\Lambda^\bullet(L)$ is an object of $\Oc$. Another class of objects
is given by the {\em (induced) Verma modules}
$
V_A = \Ind_\pen^\gen A
$
associated to finite-dimensional representations $A$ of $\men = \gl(m|n)$. 
Recall that $V_A$ is free over $U(\ol\nen)$:
\be\label{eq:Verma=free}
V_A\,\,  \= \,\,  U(\ol\nen) \otimes_\CC A \quad\text{as an $U(\ol\nen)$-module}. 
\ee
Here is our main result in its abstract (existence)  formulation.

\begin{thm}\label{thm:existence}
There is a resolution of $\gen$-modules in $\Oc$
\[
\begin{gathered}
\Vc^\bullet  = 
\bigl\{ 
\cdots \lra \Vc^{-2} \lra \Vc^{-1} \lra \Vc^0
\bigr\} 
\buildrel \on{qis}\over \to \Lambda^\bullet(L), \quad \text{where}
\\
\Vc^{-p} = V_{\Lambda^p(S^2(L))\otimes \CC_{-1/2}} =
 \bigoplus_{\alpha\in Q_{-1}(2p)\, \cap\,  \YD (m|n)} V_{\Sigma^\alpha \otimes \CC_{-1/ 2} }. 
\end{gathered} 
\]
\end{thm}

In fact, we will construct such a resolution explicitly.

\begin{exas}
(a) Let $n=0$. Then $\gen=\o(2m+1)$ and the spin module  $\Lambda^\bullet(L)$ 
if finite-dimensional with highest weight $\lambda= \bigl({1\over 2}, \cdots, {1\over 2})$. This weight is  equal to $\CC_{+1/2}$ in the notation of 
\eqref{eq:C_eps} and corresponds to the 
  {\em highest} vector; such vectors span 
 $\Lambda^n(L)$. The vector  
  $1\in\Lambda^\bullet(L)$ is, in these conventions,
   the {\em lowest} vector and  so has weight $\CC_{-1/2}$ as stated
   in Proposition \ref{prop:spin/weil}(b). 
   In this case $\Vc^\bullet$ coincides with the classical parabolic  Bernstein-Gelfand-Gelfand (BGG)
of this finite-dimensional module, see \cite{humphreys}.  An explicit comparison of the BGG
recipe $\bigoplus_w V_{w(\lambda+\rho)-\rho}$  with the form that Theorem 
\ref{thm:existence} gives in this case 
can be found in \cite[Prop.4.6]  {SS-spin}. As any classical BGG resolution, $\Vc^\bullet$ is finite in
this case, which can be seen from vanishing of $\Lambda^p(S^2(\CC^m))$ for
$p>m(m+1)/2$. 

\vskip .2cm

(b) Let $m=0$. Then $\gen = \osp(1|2n)$ and $\Lambda^\bullet(L) = S^\bullet(\CC^n)$ is infinite-dimensional, so the classical BGG theory is not applicable. In this case Theorem 
\ref{thm:existence} was proved in \cite{kapranov} by using the results of Gorelik
\cite {gorelik-typical} 
to connect $\Oc^{\gen,\pen}$ to $\Oc^{\gen_\0, \pen_\0}$,  the category $\Oc$ for $\gen_\0=\sp(2n)$ and then
the results of Enright-Shelton \cite {enright-shelton} to connect the relevant (irregular) block for 
$\Oc^{\gen_\0, \pen_\0}$ to the regular block of the category $\Oc$ for the orthogonal
group, for which the classical BGG theory applies. In this case the resolution is also
finite, as $\Vc^{-p} = V_{\Lambda^p(\Lambda^2(\CC^n))\otimes \CC_{-1/2}}$ vanishes for
$p>n(n-1)/2$. 

\vskip .2cm

(c) If $m,n>0$, then the resolution $\Vc^\bullet$ is infinite, as 
$Q_{-1}(2p) \,\cap\, \YD(m|n)\neq \emptyset$ for all $p\geq 0$. 
\end{exas}

\begin {rem}
Because of \eqref {eq:Verma=free}, the resolution $\Vc^\bullet$ gives a free resolution
of $\Lambda^\bullet(L)$ as a $U(\ol\nen)$-module. For $m=0$, 
the existence of such  free resolution
of the  shape given in Theorem \ref{thm:existence}  was established in
\cite[Prop.4.9]{SS-spin}. 
\end{rem}

\vfill\eject

\paragraph{An explicit construction of differentials: the Lie algebra chain complex.}\label{par:exp-con}
 Let  $U=U(\ol \nen)$ be the enveloping algebra of $\ol\nen$.
 As $\ol\nen = L\oplus \Lambda^2 L$
is the free nilpotent Lie superalgebra of degree $2$, the symmetric algebra $Z:= S^\bullet(\Lambda^2 L)$
lies in the center of $U$.  Thus  $Z$ parametrizes skew-symmetric forms on $L$ and so
the whole $U$ can be considered as the {\em universal Heisenberg algebra} generated by  $L$,
with $Z$ serving as the ring of coefficients. The Poincar\'e-Birkhoff-Witt theorem  for $\ol\nen$
implies that 
 we have the {\em symmetrized (Weyl)  normal form}
for elements of $U$:

\begin{prop}
(a) For each $p\geq 0$ the map  $S^p(L) \to U$ given by the symmetrized
product  is injective, so we can consider $S^p(L)$ as a subspace of $U$. Further, this extends to
an embedding $S^\bullet(L)\subset U$.

\vskip .2cm

(b) Each element of $U$ can be uniquely written as $sz=zs$ with $s\in S^\bullet(L)$ and $z\in Z$. 
This gives an identification $U\= S^\bullet(L)\otimes_\CC Z$ as a vector space. 
\qed
\end{prop}
\noindent  Consider  now the subspace  
\[
\aen \, = \,  S^2(L)\cdot Z \,  =\,  Z\cdot S^2(L) \,\subset \, U.
\]
By the above, $\aen \,\= \, S^2(L)\otimes_\CC Z$, so it is a free $Z$-module of finite rank. 

\begin{prop}\label{prop:moyal}
$\aen$ is a Lie sub(super)algebra in $U_\Lie$. The bracket in $\aen$ is $Z$-linear, so $\aen$ is
a Lie algebra over $Z$. 
\end{prop}

\noindent{\sl Proof:} The first statement  is a general property of Heisenberg algebras: the space of ``quadratic
Hamiltonians''  is closed under the Lie bracket. More precisely, direct computation
using the fact that the (super) commutators of the generators lie in the center gives 
for the (super) commutator of two (super) anticommutators of any four generators
$a,b,c,d\in L$  the answer
\[
\begin{gathered} 
\bigl[ \{a,b\}, \{c,d\}\bigr]\,\,  =\,\,  2\cdot  \biggl( (-1)^{(|b|+ |c|)\cdot |d|} \{a,d\} [b,c] +
\\
+ (-1)^{|b|\cdot |c|} \{a,c\} [b,d]  + (-1)^{(|a|+|c|) \cdot |d|} \{b,d\} [a,c]
+ (-1)^{|a|\cdot(|b|+  |c|)} \{b,c\} [a,d]\biggr). 
\end{gathered} 
\]
 This identity can be seen as a version of the the classic Moyal formula from the theory of deformation quantization, 
 saying that the
 commutator of two quadratic Hamiltonians with respect to the Moyal product is equal to their  Poisson bracket. 
  The second statement follows from this because
the coefficients at the anticommutators in this formula are  commutators
which belong to $Z$, so $\aen$ is a Lie  superalgebra over $Z$.
\qed

\vskip .2cm

Therefore the Lie superalgebra $\aen$ acts on the super vector space $U$ on the left by left multiplication in $U$
and on the right by right multiplication in $U$. It will be convenient for us to use the right action. 

\vskip .2cm

Let $(\FMod_Z, \otimes_Z)$ be the category of free $Z$-(super) modules,
i.e., of $\ZZ/2$-graded free $Z$-modules,
equipped with the symmetric monoidal structure
$\otimes_Z$. Here the symmetry data for $\otimes_Z$ are given by the usual Koszul sign rule for super-modules.
In particular, for an object $M\in \FMod_Z$ we can speak of its symmetric and exterior powers over $Z$,
denoted $S^p_Z(M)$ and $\Lambda^p_Z(M)$. 

\vskip .2cm

Note that both $U, \aen\in\FMod_Z$ and the right action
 can be written as a morphism
$\act: U\otimes_Z \aen \to U$ in  $\FMod_Z$ given by  $u\otimes a \mapsto ua$.

\vskip .2cm

Let $\Cc^\bullet$ be the Chevalley-Eilenberg chain complex ({\em in the category $\FMod_Z$}) calculating the Lie algebra homology of $\aen$ with coefficients in $U$,
graded so that it sits in negative degrees and the differential 
$d= d^\Lie$ has degree $+1$:
\be
\begin{gathered}
\Cc^\bullet  = C^\Lie_{-\bullet}(\aen, U)_Z \,=\, 
\bigl\{ \cdots  \buildrel d_3\over \lra U\otimes_Z \Lambda^2_Z \aen 
\buildrel d_2\over \lra U\otimes_Z \aen\buildrel d_1\over  \lra U\bigr\} 
\\
=\hskip 2.7cm \bigl\{ \cdots \to U\otimes_\CC \Lambda^2(S^2L) \to U\otimes_\CC S^2L \to U\bigr\}. 
\end{gathered}
\ee
 Thus  \[
 d_1=\act, \quad d_2\bigl((u\otimes (a_1\wedge a_2)\bigr) =  (ua_1)\otimes a_2 - (-1)^{|a_1|\cdot |a_2|}
 (ua_2)\otimes a_1  - u\otimes [a_1, a_2]
 \]
 and so on, see \cite{fuchs, garcia}.  Our conventions imply that  $\Cc^\bullet$ is a complex of free {\em  left}  $U$-modules,
 as the differentials, coming from the {\em right} $\aen$-action, do not affect the left side of the left factors of
 the tensor
 products. In particular, each $H^{-i}(\Cc^\bullet) = H^\Lie_i(\aen, U)$ is a left $U$-module.

 Now, Theorem \ref{thm:existence} will follow from the next two more precise results which
 will be proved in the rest of the paper.
 
 \begin{thm} \label {thm:exactness}
 $\Cc^\bullet$ is a left resolution of $\Lambda^\bullet(L)$ as a module over $U = U(\ol\nen)\subset U(\gen)$, with the augmentation
  morphism $U\to\Lambda^\bullet(L)$ being the natural algebra
 homomorphism obtained by factorizing $U$ by the two-sided ideal generated by
 $S^2(L)$.  That is,
 \[
H_i^\Lie(\aen, U) = \begin{cases} \Lambda^\bullet(L), &\text{ if } i=0;
\\
0, & \text{ if } i\neq 0. 
\end{cases} 
 \]
 \end{thm}
 
 Note that $\Cc^\bullet$ is a $\GL(L)$-equivariant complex and the differentials are $\GL(L)$-invariant. 
 This makes $\Cc^\bullet$ a complex of $\ol\pen$-modules. 
  Let us  modify the $\ol\pen$-module structure
 by multiplying $\Cc^\bullet$ by the $1$-dimensional representation 
 $\CC_{-1/2}$ of $\men = \gl(L)$. By \eqref{eq:Verma=free} we have an identification
 \[
 \Cc^{-p}\otimes \CC_{-1/2}\= V_{\Lambda^p(S^2(L))\otimes\CC_{-1/2}}. 
 \]
 of $\ol\pen$-modules with  the RHS of this identification being  in fact a $\gen$-module
 which makes  each $\Cc^{-p}\otimes \CC_{-1/2}$ into a $\gen$-module as well.   
 As a vector space, it is identified with $\Cc^{-p}$. 
 
 \begin{thm}\label{thm:linearity}
With respect to these  $\gen$-module structures, each differential
\[
d_p\otimes \Id: \Cc^{-p}\otimes\CC_{-1/2} \to \Cc^{-p+1}\otimes\CC_{-1/2}
\]
 as well
as the augmentation map $\Cc^0\otimes \CC_{-1/2}\to \Lambda^\bullet(L)$
are morphisms of $\gen$-modules, i.e., they are  $U(\gen)$-linear. 
Therefore 
 \[
 \Vc^\bullet = \bigl(\Cc^\bullet\otimes \CC_{-1/2},
 d^\Lie\otimes\Id\bigr)
 \]
 is a resolution claimed in  Theorem 
 \ref {thm:existence}.
 \end{thm}

 \vfill\eject
 
 \section {Proof of exactness } 
 
 Here we prove Theorem \ref{thm:exactness}. We start with:

 \paragraph{The rightmost term.} 
 By construction,
 \[
 H^0(\Cc^\bullet) = H_0^\Lie(\aen, U) =  \Coker \{\act: U\otimes_Z \aen \to U\} = U/U\aen = U\bigl/(U\cdot S^2(L)).
 \]
 
 \begin{prop}\label{prop:last-term}
 (a) The left ideal $U\cdot S^2(L) $  in $U$ coincides with the right ideal $S^2(L)\cdot U$
 and so is  equal to the $2$-sided ideal 
 $U\cdot S^2(L)\cdot U = \Ker \{U\to \Lambda^\bullet(L)\}$. 
 
 \vskip .2cm
 
 (b) Therefore 
 the space $H^0(\Cc^\bullet)$ is
 isomorphic to $\Lambda^\bullet(L)$ as a left $U$-module. 
 \end{prop}
 
 \noindent{\sl Proof:} (a) It is enough to show that 
 $L\cdot S^2(L) = S^2(L)\cdot L$
 as subspaces in $U$.  By the PBW theorem $U$ splits, as a $\GL(L)$-module, into the
 sum
 \[
 U \, \= \, S^\bullet(L\oplus \Lambda^2(L))\,\,  = 
 \,\,
 S^\bullet(L) \otimes S^\bullet(\Lambda^2(L)) \= 
  \bigoplus_{\alpha\in \YD(m|n)} \Sigma^\alpha(L).
 \]
 That is, each nonzero $\Sigma^\alpha (L)$ enters into
 $U$ exactly once. 
 See, e.g., \cite[\S4.3]{SS-spin} for the case of purely even $L$; the case of arbitrary $L$
 follows from the resulting identification of polynomial functors on superspaces. 
 Consider the $\ZZ$-grading 
  $U = \bigoplus_{d\geq 0} U^d$ on $U$
  obtained by putting $L$ in degree $1$
 and $\Lambda^2(L)$ in degree $2$. Then 
 \[
 L\cdot  S^2(L), \,  S^2(L)\cdot L \quad \subset \quad  U^3 \,\, \= \,\, S^3(L) \oplus \Sigma^{2,1}(L) 
 \oplus\Lambda^3(L)
 \]
 are $\GL(L)$-submodules and so each of them is the sum of some subset of the summands
 on the right. Since each of  these subspaces is  the image of $L\otimes S^2(L)$
 under a $\GL(L)$-equivariant map and $L\otimes  S^2(L)$ does not contain $\Lambda^3(L)$, 
 they must be sum of some subsets  of 
 $\{ \Sigma^{2,1}(L), S^3(L)\}$. We claim that these
 subsets are the same (in fact, this common subset is the full set of two summands).
 Indeed, as any enveloping algebra, $U=U(\ol\nen)$ carries the antipode, a superalgebra anti-involution
 $\sigma$ defined by $\sigma(x)=-x$, $x\in\ol\nen$. By construction, this anti-involution
 is $\GL(L)$-invariant, so it preserves each of the summands $\Sigma^\alpha(L)\subset U$
 and  is scalar 
 (in fact, $\pm 1$) on this summand. But since $\sigma$ is an anti-automorphism,
 it takes $L\cdot S^2(L)$ to $S^2(L)\cdot L$ and so these two subspaces coincide. 
 This proves (a).  Part (b) follows. \qed

  \paragraph{$U$ and $\Cc^\bullet$ as polynomial functors of $L$: reduction to $L$ odd.}\label{par:U-poly}
 To emphasize the dependence of our constructions on $L=\CC^{m|n}$ which
 can be any finite-dimensional superspace, we write $\ol\nen$ as
 $\ol\nen(L) = L\oplus S^2(L)$. This makes $\ol\nen$ into  a  polynomial functor on
 $\SVect^\fd$, i.e., an 
 object of $\SPF$. Therefore (as noted in \cite{getzler} in the purely even case)
  the correspondence
 $L\to U(L) :=  U(\ol \nen(L))$ is also an  object of $\SPF$.  Indeed, $U(\ol\nen(L))$
 is obtained from $\ol\nen(L)$ by tensor operations (quotienting the tensor algebra
 by the standard relations)
  which have their counterparts
 in the category $S\Mod$ with its  tensor structure $*$. 
 
 
 For the same reason
  the complex $\Cc^\bullet =\Cc^\bullet(L) := C_{-\bullet}^\Lie(\aen,  U)$ is, in its dependence on $L$,  also a complex
 of polynomial functors.  In fact, all the polynomial functors above are even. In addition,
 the augmentation  morphism $U = \Cc^0 \to \Lambda^\bullet(L)$ is also a morphism of  polynomial
 functors. Therefore 
 Proposition \ref {prop:exactness-even} reduces the proof of Theorem \ref{thm:exactness}
 to the case when $L = \CC^{0|n} $ is purely odd which we assume from now on. 
 
 \vskip .2cm
 
 It is convenient to write 
 $L = \Pi(L')$ where
 $L'=\CC^n$ is purely even. Thus the Lie superalgebra $\ol\nen$
 has the form $L' \oplus  S^2(L')$
 where $L'$ sits in degree $\1$. The super-space 
 $\Lambda^\bullet(L) = S^\bullet(L')$ is the Weil representation of
 $\osp(1|2n)$, see Example \ref{exas:spiun/osc}(b). 
  Further,  $S^2(L)=\Lambda^2(L')$ is  realized inside $U$ as the
  span of expressions $ab-ba$, $a,b\in L'$, which are, in this case,
 super-anticommutators \eqref{eq:super-anticom}.

 \paragraph{Localization on the spectrum of the center of $U$.} 
 Let $X:=\Spec (Z)$ be the spectrum of the center of $U$.
 For any $Z$-module $M$ we denote by $\ul M$ the corresponding
 quasi-coherent sheaf on $X$. 
 
 \vskip .2cm

 As $Z= S^\bullet(\Lambda^2(L)) = S^\bullet(S^2(L'))$,
 we view   $X$ as $S^2(L')^*$, the space of quadratic forms
 on $L'$. The algebra $U$, being the universal Clifford algebra,
 is a free module over $Z$ of rank $2^n$, so $\ul U$
 is a vector  bundle on  $X$ of this rank. 
 
 \vskip .2cm
 
 Further, each term   of the complex
 $\Cc^\bullet$ is also free over $Z$ of finite rank, namely
 $\Cc^{-p}$ is free of rank $2^n\cdot {n(n-1)/2\choose p}$. 
 Thus each $H^i(\Cc^\bullet)$ is a finitely generated $Z$-module. 
  Proposition \ref{prop:last-term} gives
 that $H^0(\Cc^\bullet)=S^\bullet(L')$, and we need to prove that
 other $H^i(\Cc^\bullet)$ vanish. 
 
 \vskip .2cm
 
 Let $X_{\leq r}\subset X$ be the locus of quadratic forms of rank $\leq r$. 
 We will also write $X_{=r}=X_{\leq r} \- X_{\leq r-1}$,
 $X_{>r} = X\- X_{\leq r}$ etc. Our strategy will be to show that the
 morphism $\ul\Cc^\bullet \to \ul{S^\bullet(L')}$
 of complexes of coherent sheaves on $X$, 
 coming from Proposition \ref {prop:last-term},
 is a quasi-isomorphism on the   open sets $X_{\geq r}$ for  successively
  smaller $r$. 
 
 \vskip .2cm
 
  Let 
 $\pi: (L')^* \to S^2(L')^*$
 is the map $y\mapsto y^2$. In the matrix form, viewing $y$ as a column vector
 and $S^2(L')^*$ as the space of symmetric matrices,
 we can write  $\pi(y) = y\cdot y^t$. 
 Thus the image of $\pi$ is $X_{\leq 1}$ with $\pi$ being $2:1$ over $X_{=1}$
 and $1:1$ over $0 = X_{=0}$. 
 The following is straightforward from the definitions.
 
 \begin{prop}\label{prop:ulS(L')}
 The coherent sheaf  $\ul{S^\bullet(L')}$ on $X = S^2(L')^*$ associated to
 $S^\bullet(L')$ as a $Z$-module, has the form
  $\ul{S^\bullet(L')} = \pi_* \Oc_{(L')^*}$. \qed
 \end{prop}

 \paragraph{ Exactness over rank $\geq 2$}
  Proposition \ref {prop:ulS(L')} implies, in
 particular, that  the $Z$-module $S^\bullet(L')$ is supported on $X_{\leq 1}$.
 
 \begin{prop}\label{prop:ex-rk-2}
 All the cohomology modules $H^i(\Cc^\bullet)$ are supported on $X_{\leq 1}$. 
 \end{prop}
 
 \noindent{\sl Proof:} As $\ul\Cc^\bullet$ is a complex of vector bundles
 on $X$ of finite rank, it suffices to prove that for any $\CC$-point
 $q\in S^2(L')^*$, i.e., for any quadratic form on $L'$, which has
  rank $\geq 2$,
 the corresponding  complex of fibers at $q$, i.e., 
 \[
 \ul \Cc^\bullet_q \,=\, \Cc^\bullet\otimes_Z \CC_q \,=\, C_{-\bullet}^\Lie
 (\aen_q, \Cliff(q)),
 \]
 is exact. Here:
 \begin{itemize} 
 
 \item $\CC_q$ is the $1$-dimensional $Z$-module corresponding
 to $q$; 
 
 \item  $\Cliff(q)= \Cliff(L', q) $ is the Clifford algebra of the quadratic form $q$ 
 on $L'$;
 
 \item  $\aen_q\subset \Cliff(q)$ is Lie algebra of quadratic Clifford
 operators spanned by the commutators $ab-ba$, $a,b\in L'$
 and identifed with $\Lambda^2(L')$ as a vector space. 
 
 \end{itemize}
 
 \noindent 
 Further, it is enough to prove this for the case when the rank of $q$
 is exactly $2$. Indeed, the complex $\Cc^\bullet$
 and therefore each $H^i(\Cc^\bullet)$ being $\GL(L')$-equivariant
 $Z$-modules, each 
 $\on{Supp}  H^i(\Cc^\bullet)$
 is a closed union of $\GL(L')$-orbits in $X$. These orbits
 are   the $X_{=r}$, with $X_{\leq r}$ being the closure of $X_{=r}$.
  So if the support
 does not meet $X_{=2}$, it does meet any larger orbits. 
 
 \vskip .2cm
 
 So we assume that $q$ has rank $2$. By $GL(L')$-equivariance,
 it is enough to take $q(x_1, \cdots, x_n) = x_{n-1}^2 + x_n^2$. 
 Thus $\Cliff(q)$ is generated by the  elements $\xi_1, \cdots, \xi_{n-2},
 e_1, e_2$ which anticommute with each other and satisfy $e_i^2=1$,
 $\xi_j^2=0$. In other words, it can be written as 
 $\Cliff_2[\xi_1,\cdots, \xi_{n-2}]$ where $\Cliff_2$ is generated by
 $e_1, e_2$. The Lie algebra $\aen_q$ is spanned by 
 \[
 e_1e_2 = {1\over 2} [e_1, e_2], \quad e_i\xi_j = {1\over 2} [e_i, \xi_j], 
 \quad \xi_i\xi_k = {1\over 2} [\xi_i, \xi_k]. 
 \]
 The element $e_1e_2$ generates the orthogonal Lie algebra $\oen(2)$
 which is abelian and identified with $\CC$. 
 
 \vskip .2cm 
 
 In $\Cliff(q)$ we have  $(e_1e_2)^2=-1$. Therefore 
  the eigenvalues of  the action of $e_1e_2$ on $\Cliff(q)$ by right multiplication are $\pm i$. This implies that the adjoint action
 of $e_1e_2$ in $\End_{\mathbb C}(\Cliff(q))$ has eigenvalues $\pm 2i,0$. In particular, $\operatorname{ad}_{e_1e_2}$ on $\aen\subset \End_{\mathbb C}(\Cliff(q))$ has eigenvalues $\pm 2i,0$. We obtain that the action of  $e_1e_2$ on each term of the Chevalley--Eilenberg complex $\Cc^p=\Lambda^p(\aen)\otimes \Cliff(q)$ has eigenvalues in the set $(\{\pm 1\}+2\mathbb Z)i$ 
 which does not contain $0$. Since the action of $\aen$ and, in particular, of $e_1e_2$ on the cohomology must be trivial, we obtain $H_{\bullet}^\Lie
 (\aen_q, \Cliff(q))=0$. 
  Proposition \ref {prop:ex-rk-2} is proved. 
 
 \paragraph{Rank $\leq 1$: numerics.} Theorem  \ref{thm:exactness},
 if true, would imply the identification
 \[
 H^{- \bullet} (\ul \Cc^\bullet_q)  \= \Tor_\bullet^Z (\CC_q, S^\bullet(L'))
 \]
 for any quadratic form $q$ on $L'$. This is because $\Cc^\bullet$ is
 a complex of free $Z$-modules and so $\ul\Cc_q$,  the complex of fibers at $q$, 
 computes the hyper-Tor between $\CC_q$ and $\Cc^\bullet$. 
So as a step towards the full proof we establish the following numerical fact.

\begin{prop}\label{prop:numerics-rk1}
If $q$ has rank $\leq 1$, then for any $i>0$ we have 
\[
\dim_\CC H^{-i} (\ul \Cc^\bullet_q) \, = \, \dim_\CC  \Tor_i^Z (\CC_q, S^\bullet(L')). 
\]
\end{prop}
 
 \noindent{\sl Proof:} Consider first the case when $q$ has rank $1$,
 so $q(x_1,\cdots, x_n) = x_n^2$. Thus 
 $\Cliff(q) = \Cliff_1 [\xi_1,\cdots,\xi_{n-1}]$ , where $\Cliff_1  = \CC[e]/(e^2-1)$. 
 The Lie algebra $\aen_q$ is spanned by the  $e\xi_i$ and $\xi_i\xi_j$ with
 the   only nonzero bracket being
 \[
 [e\xi_i, e\xi_j] = -2 \xi_i \xi_j
 \]
 Denoting by $\CC^{n-1}$ the space spanned by the $e\xi_i$,
 the space spanned by the $\xi_i \xi_j$ can be written as $\Lambda^2\CC^{n-1}$,
 and $\aen_q = \CC^{n-1} \oplus \Lambda^2 \CC^{n-1}$
 is the free $2$-step nilpotent algebra on $\CC^{n-1}$, i.e., the
 nilpotent radical in the  parabolic subalgebra  of $\oen(2n-1)$
 preserving an isotropic $\CC^{n-1}$. In the notation of
 \S \ref {par:U-poly}, we can say that $\aen_q = \ol\nen(\CC^{0|n-1})$.

  With respect to this identification, we further see that $\Cliff(q)$ is, as an
  $\aen_q$-module,  isomorphic to the direct sum of two copies of
  $\Lambda^\bullet(\CC^{n-1})$ as a module over $\CC^{n-1} \oplus \Lambda^2 \CC^{n-1}$. Indeed,  let us decompose $\Cliff(q)$,
  of dimension $2^n$,  into the sum of two eigenspaces of the operator of
  left multiplication with $e$, the eigenvalues being $\pm 1$.
  These eigenspaces are $\aen_q$-submodules of dimension $2^{n-1}$,
  each isomorphic to $\Lambda^\bullet(\CC^{n-1})$ as above. Therefore
  \be\label{eq:H^iC_q}
  \dim H^{-i}(\ul \Cc^\bullet_q) = 2 \dim H_i^\Lie (\CC^{n-1} \oplus \Lambda^2\CC^{n-1}, \,\Lambda^\bullet \CC^{n-1}) = 2 {n(n-1)/2 \choose i}. 
  \ee
 Here the latter equality comes from the identification of the middle
 term given in \cite[Prop.4.6]{SS-spin}. 
 
 \vskip .2cm
 
 On the other hand, Proposition \ref {prop:ulS(L')} and the fact that $\pi$ is
 $2:1$ over $X_{=1}$ imply that the sheaf  $\ul{S^2(L')}$ is,
   Zariski locally near $q$,
  isomorphic to the direct sum of
 two copies of $\Oc_{X_{=1}}$. Therefore
 \[
 \dim \Tor_i^Z (\CC_q, {S^2(L')}) = 2\dim \Tor_i^Z(\CC_q, \Oc_{X_{=1}})
 = 2 { \codim_X X_{=1} \choose i}  =  2 {n(n-1)/2 \choose i},
 \]
 as it follows from the Koszul resolution of the structure sheaf of the
 smooth subvariety $X_{=1}$ near $q$. This establishes Proposition
 \ref {prop:numerics-rk1} for $q$ of rank $1$.
 
 \vskip .2cm
 
 We now consider the case  $q=0$. Then $\Cliff(0) = \Lambda^\bullet(L')$
 and $\aen_0$ is the vector space $\Lambda^2(L')$ with trivial bracket. 
 Thus
 \[
 H^{-i} (\ul \Cc^\bullet_0) = \Tor_i^{S^\bullet(\Lambda^2L')} (\CC_0, 
 \Lambda^\bullet  L' ), 
  \]
 while the RHS of the formula in Proposition
 \ref {prop:numerics-rk1} involves 
 similar data but 
 with the exterior powers instead of 
  the symmetric ones. To identify their dimensions, we recall
  the result of  J\'ozeffiak-Pragacz-Weyman \cite[Th.3.19]{weyman-resol},
  see also \cite[Th.6.3.1]{weyman-book},
 that gives an identification of $\GL(L')$-modules
  \be\label{eq:weyman-form}
  \Tor_i^{S^\bullet(S^2L')} (\CC_0, S^\bullet L') \, \= \,
  \bigoplus_{\alpha = \alpha ^\vee, \atop (|\alpha| - s(\alpha))/2=i} \Sigma^\alpha L'.
  \ee
  Here the summation is over self-dual (symmetric) Young diagrams and 
  $s(\alpha)$ is the size of the Durfee square of $\alpha$, i.e., of the
  maximal sub-square inscribed in $\alpha$. Now, the LHS of 
  \eqref {eq:weyman-form} is a polynomial functor which can be applied
  to any vector superspace $L'$. In particular, it can be applied to 
  $L=\Pi (L')$ and the identification   \eqref {eq:weyman-form}
  will still hold if we understand  $\Sigma^\alpha (\Pi(L'))$ in the
  super-sense, as $\Sigma^{\alpha^\vee}(L')$. 
  But the diagrams in  \eqref {eq:weyman-form} are self-dual,  and so
  we get the identical formula
  \[
H^{-i}(\ul\Cc^\bullet_0)
  \, = \,
  \Tor_i^{S^\bullet(\Lambda ^2L')} (\CC_0, \Lambda^\bullet L') \,\,\,\, \= \,
  \bigoplus_{\alpha = \alpha ^\vee, \atop (|\alpha| - s(\alpha))/2=i} \Sigma^\alpha L' 
  \,\,\,\, 
  \=  \,\,\,\,   \Tor_i^{S^\bullet(S^2L')} (\CC_0, S^\bullet L'),
     \] 
 thus identifying not only the dimensions but also the $\GL(L')$-module
 structures. This finishes the proof of Proposition  \ref {prop:numerics-rk1}. 
 
 \paragraph{Rank $\leq 1$: exactness.} We now finish the proof of 
 Theorem \ref{thm:exactness}, i.e.,  prove  that $H^i(\Cc^\bullet)=0$
 for all $i<0$. By Proposition \ref {prop:ex-rk-2} these $H^i$ are supported
 on $X_{\leq 1}$ anyway. So we eliminate them in two steps.
 
 \vskip .2cm
  First we
 prove that the $H^i(\Cc^\bullet)$, $i<0$, if any, are supported at $0$.
 Suppose the contrary, i.e., that there is $m>0$ such that 
  $\on{ supp} H^{-m}(\Cc^\bullet) = X_{\leq 1}$ (by $\GL(L')$-equivariance,
  the support must be either $0$ or  the entire $X_{\leq 1}$). 
  Let $m$ be the largest integer with this property.    
  Denote by $N=n(n-1)/2$ the codimension of $X_{\leq 1}$ in $X$
  and let $q\in X_{=1}$ be a quadratic form of rank $1$.

  \begin{lem}\label{lem:box}
  If  $j<0$ is such that $H^j(\Cc^\bullet)$ is supported on $X_{\leq 1}$,
  then:
  
  \vskip .2cm
  
  (a) $\Tor_i^Z(\CC_q, H^j(\Cc^\bullet)) =0$ for $i>N$.
  
  \vskip .2cm
  
  (b) $\Tor_N^Z(\CC_q, H^j(\Cc^\bullet)) \neq 0$
  \end{lem}
  
  \noindent {\sl Proof:}   
  Let $I\subset \Oc_X$ be the ideal of $X_{\leq 1}$. Any $H^{i}(\ul \Cc^\bullet)$,
  being a coherent sheaf on $X$ supported on $X_{\leq 1}$, 
  has its $I$-adic filtration 
  $\{I^k  H^{i}(\ul \Cc^\bullet)\}$ terminating. The quotients of this filtration
  are coherent sheaves on  $X_{\leq 1}$. 
  Restricted to $X_{=1}$, they  are, by  $\GL(L')$-equivariance, vector bundles
  on $X_{=1}$ so locally isomophic to direct sums of several copies of
  $\Oc_{X_{=1}}$. Now, $ \Tor_i^{\Oc_{X_{\geq 1}}}(\CC_q, \Oc_{X_{=1}})$
  is calculated by a Koszul resolution of $\Oc_{X_{=1}}$ near $q$
  and its dimension is ${N\choose i}$. This proves (a).
  
  \vskip .2cm
  
  To prove (b), we note that the deepest nonzero term of the $I$-adic filtration
  is, by the above, a vector bundle on $X_{=1}$ near $q$. Thus we have
  an embedding $\Oc_{X_{=1}} \hra H^j(\ul\Cc^\bullet)$ near $q$
  which by (a) and by the long exact sequence of $\Tor$'s, induces an embedding 
  $\Tor_N^Z (\CC_q, \Oc_{X_{=1}})\hra \Tor_N^Z (\CC_q, H^j(\ul \Cc^\bullet))$.
  Now, the source of this embedding has dimension $1$ by the Koszul
  resolution calculation above, so the target is nonzero.  \qed. 
  
  \vskip .2cm
  
  Consider now the spectral sequence
  \[
  E_2^{ij} = \Tor_{-j}^Z (\CC_q, H^i(\Cc^\bullet)) \,\,\Rightarrow \,\,
  \Tor_{-j+i}^Z (\CC_q, \Cc^\bullet) \,=\, H^{j-i}(\ul\Cc_q), \quad i,j\leq 0, 
  \]
  see Fig. \ref {fig:SS}. Its differential $d_r$, $r\geq 2$, has
 bidegree $(-r+1, r)$.

 \begin{figure}[h] 
   \centering
   \btp[scale=0.3]
 
 \draw (-2,0) -- (33,0); 
  \draw (-2,3) -- (33,3); 
   \draw (-2,8) -- (33,8); 
    \draw (-2,11) -- (33,11); 
     \draw (-2, 14) -- (33,14); 
     
   \draw (0, -2) -- (0, 14);   
    \draw (11, -2) -- (11, 14);   
       \draw (22, -2) -- (22, 14);   
          \draw (33, -2) -- (33, 14);   
          
  \node at (5.5, 1.5){\small$\Tor_N(\CC_q, H^{-m}\Cc^\bullet)$};   
  \node at (16.5, 1.5){$\cdots\cdots$}; 
  \node at (27.5, 1.5){$\Tor_N(\CC_q, H^0\Cc^\bullet)$}; 
  
   \node at (5.5, 5.5){$\cdots\cdots $};   
    \node at (5.5, 9.5){\small$\Tor_1(\CC_q, H^{-m}\Cc^\bullet)$};   
     \node at (5.5, 12.5){\small$\Tor_0(\CC_q, H^{-m}\Cc^\bullet)$};   
     
        \node at (16.5, 5.5){$\cdots\cdots $};   
           \node at (16.5, 9.5){$\cdots\cdots $};   
   \node at (16.5, 12.5){$\cdots\cdots $};   
   
      \node at (27.5, 5.5){$\cdots\cdots $};   
  \node at (27.5, 9.5){\small$\Tor_1(\CC_q, H^{0}\Cc^\bullet)$}; 
    \node at (27.5, 12.5){\small$\Tor_0(\CC_q, H^{0}\Cc^\bullet)$}; 
    
    \draw[fill, color={gray}, opacity =0.2] (0,0) -- (11,0) -- (11,3) -- (0,3) -- (0,0); 
    
    \draw [->, line width =0.8] (0.5, 2.5  ) -- (-6,8); 
    \node at (-4,5){$d_2$ } ; 
 
 \etp
 \caption{The lower left term survives in $E^\oo$. }
 \label{fig:SS}
 \end{figure}
 
  By Lemma \ref {lem:box}, nonzero terms  $E^{ij}_2$
 are situated in the  box $-m\leq i \leq 0$, $-N\leq j \leq 0$,
 and the term $E_2^{-m, -N}$ in the lower left corner is nonzero. 
 We now note that the differentials $d_2, d_3$ etc. applied to this term, 
 have their target outside the box of nonzero terms and similarly
 such differentials targeting it, must come from outside the box.
 This implies that $E_2^{-m, -N}$ survives in $E_\oo$ and contributes
 to $\Tor_{m+N}^Z(\CC_q, \Cc^\bullet)$, i.e., to $H^{-m-N}(\ul\Cc_q)$,
 the $(-m-N)$th cohomology of the complex of fibers at $q$.
 But we know by \eqref {eq:H^iC_q} that $H^{j}(\ul\Cc_q)=0$
 for $j<-N$. Thus $m>0$ is impossible and all $H^j(\Cc^\bullet)$,
 $j<0$, if any, must be supported at $0$. 
 
 \vskip .2cm
 
 To eliminate the possibility of support at $0$, we proceed in a similar
 way to the above. Let $D=n(n+1)/2 = \dim(X)$ and $J\subset Z$
 be the maximal ideal of $0$. Then 
 $\Tor_i^Z(\CC_0, \CC_0)$ has dimension ${D\choose i}$. 
 Any nonzero $H^{j} (\Cc^\bullet)$, $j<0$, being supported at $0$,
 has the $J$-adic filtration finite with each quotients being a  direct sum of
 several copies of $\CC_0$. This implies that 
 $\Tor_i^Z(\CC_0, H^{j} (\Cc^\bullet))= 0$ for $i>D$. Similarly, looking at the
 deepest nonzero layer of the filtration, we find that 
 and $\Tor^Z_D(\CC_0, H^{j} (\Cc^\bullet))\neq 0$. Denoting by
 $m$ the largest integer such that $H^{-m}(\Cc^\bullet)\neq 0$,
 we consider the spectral sequence
  \[
  E_2^{ij} = \Tor_{-j}^Z (\CC_0, H^i(\Cc^\bullet)) \,\,\Rightarrow \,\,
  \Tor_{-j+i}^Z (\CC_0, \Cc^\bullet) \,=\, H^{j-i}(\ul\Cc_0), \quad i,j\leq 0, 
  \]
 with nonzero terms being restricted to the box $-m\leq i\leq 0$,
 $-D\leq j \leq 0$ with $E_2^{-m, -D}\neq 0$. By the same argument
 as above, this term will survive in $E_\oo$ and contribute to 
  $\Tor_{m+N}^Z(\CC_0, \Cc^\bullet)$, i.e., to $H^{-m-N}(\ul\Cc_0)$.
  But $H^{j}(\ul\Cc_0)=0$ for $j<-D$ by Proposition \ref {prop:numerics-rk1}.
  
  \vskip .2cm
  
  This finishes the proof of Theorem  \ref{thm:exactness}. 
  
   \vfill\eject
   
   \section {$\gen$-linearity of the differentials}
   
   Here we prove Theorem \ref{thm:linearity}.
   
   \paragraph {$\ssl(2)$-style relations in $\gen$.} 
    We retain the notation  of \eqref{eq:g-notation} as to $\gen=\osp(2m+1|2n)$,  $\men=\gl(m|n)=\gl(L)$,
    $L=\CC^{m|n}$  etc. 
     As before, we also abbreviate $U(\ol\nen)$ to $U$.
      Since  $\ol\nen = L\oplus \Lambda^2(L)$ and
    $\nen = L^* \oplus \Lambda^2 L^*$, we regard any $y\in L$ as an element of $\ol\nen$
    and any $z\in L^*$ as an element of $\nen$.

    We will use the following relations in $\gen$ which are matrix versions
    of the familiar relations in $\ssl(2)$ to which they reduce for $L=\CC$.
    Here $h\in \men$ while $y,y'\in L$ and $z\in L^*$:
    \be\label{eq:rels-in-g}
    \begin{gathered}
    [h,y] = h(y), \,\, [h,z] = -h^t (z), \quad  [y,z] = 2 (-1)^{|y|\cdot |z|} e_{zy}, \text{ where}
    \\
  e_{zy} \in \men\,\, \text{ is defined by } \,\,  e_{zy}(y') = 
   (-1)^{|y|\cdot |y'|} (z, y')\cdot y.       \end{gathered}
    \ee
   Thus $e_{zy}$ corresponds to $z\otimes y\in L^*\otimes L$ under the identification of the latter space with $\men$. 
    Note that   $\str (e_{zy}) = (z,y)$.

   \paragraph{Criterion of $\gen$-linearity.}   Let $A$ be a finite-dimensional representation of
   $\men=\gl(m|n)$ and $V$ any $\gen$-module. By definition
   \be\label{eq:verma-uni}
   \Hom_\gen (V_A, V) \= \Hom_\men (A, V^\nen). 
   \ee
 This implies that the augmentation map in Theorem \ref{thm:linearity} is $\gen$-linear, since $1\in\Lambda^\bullet(L)$
 is annihilated by $\nen$ and is an $\men$-eigenvector with character $\CC_{1/2}$. To treat the other differentials
 in $\Cc^\bullet\otimes\CC_{1/2}$, we start with the following.
 
 \begin{prop}\label{prop:verma-lin}
 Let $f: V_A\to V$ be a $\ol\pen$-linear map, i.e., an $\men$-equivariant $\ol\nen$-linear map. 
 Suppose $f(A)\subset V^\nen$. Then $f$ is $\gen$-linear. 
 \end{prop}
 
 \noindent {\sl Proof:} As $f: A\to V^\nen$ is $\men$-linear, \eqref{eq:verma-uni}  gives  a $\gen$-linear map
 $f': V_A \to V$ defined by $f'|_A = f|_A$. We claim that $f'=f$ and so $f$ is $\gen$-linear. Indeed, $f'=f$ on 
 $A\subset V_A$. Further, as a $U$-module $V_A = U\otimes_\CC A$ is free,
 so $f, f': U\otimes_\CC A\to V$ are two morphisms from a free $U$-module to another
 $U$-module which 
 coincide on the space of generators. So $f=f'$. \qed

 \paragraph{Polynomial functor reformulation. Reduction to $L$ even.} Let $A$ be
  a finite-dimensional representation of $\men$. Since the identification 
  $V_A \= U\otimes_\CC A$ is induced by the $U$-action, we denote the generators by
  $\ket{a}\in V_A,$ to correspond to   $1\otimes a$, $a\in A$ and write  $u\ket{a}\in V_A$ to correspond to $u\otimes a$, $u\in U$ under this identification. 
  
  \vskip .2cm 
  
  Consider  the $\ZZ$-grading
  $V_A=\bigoplus_{d\geq 0} V_A^d$ with $V_A^d$ corresponding to $U_A^d\otimes_\CC A$
 and $U=\bigoplus_{d\geq 0} U^d$ being the grading such that $U^1=L$.  Thus
 $V_A^d$ is spanned by the $y_1\cdots y_d \ket{a}$ for all $y_i\in L$ and $a\in A$. 
 The action of  any
 $y\in L$ increases the grading by $1$, the action of any $z\in L^*$ decreases it by $1$
 and the action of any $g\in \men$ preserves it. 
 
   \vskip .2cm 
   
   Recall that $U = U(\ol\nen(L))$ as well any graded component $U^d$
   can be considered as  a (super) polynomial functor in $L$. 
    Suppose now  that $A$  also comes from a  polynomial functor (applied to $L$
    to get a representation of $\GL(L)$
   and then differentiated to get a representation of the Lie algebra). Then
   each $V_A^d \= U^d \otimes_\CC A$ is  again a polynomial functor in $L$. 
   For example, $L\mapsto V^d_{\Lambda^p(S^2L)}$ is a polynomial functor. 
   
    \vskip .2cm
 
 By Proposition \ref {prop:verma-lin},  to prove that 
 \[
 d_p\otimes \Id: \Cc^{-p}\otimes\CC_{-1/2}  = U\otimes_\CC \Lambda^p(S^2 L)\otimes \CC_{-1/2} \lra
 U\otimes_\CC \Lambda^{p-1}(S^2 L)\otimes \CC_{-1/2} = 
 \Cc^{-p+1}\otimes \CC_{-1/2}
 \]
 is $\gen$-linear, it suffices to show that the subspace $(d_p\otimes \Id) (1\otimes \Lambda^p(S^2 L ) \otimes 1)$
 is annihilated by $\nen = L^*\oplus \Lambda^2(L^*)$. Since $\nen$ is generated by $L^*$,
 it is enough to show that this subspace is annihilated by any $z\in L^*$. 
 After untwisting by $\CC_{-1/2}$ both in source and target, this gives a 
 $GL(L)$-equivariant
 map
 \[
  \gamma_L:   L^* \otimes \Lambda^p(S^2L) \lra V_{\Lambda^{p-1}(S^2L)}^1 =  L\otimes \Lambda^{p-1}
 (S^2L), \quad z\otimes a \mapsto z d\ket {a}
 \]
 and we need to prove that $\gamma_L=0$. Consider the partial transpose
 which is a $GL(L)$-equivariant map
 \[
 \gamma_L^t : \Lambda^p(S^2L) \lra L \otimes L \otimes \Lambda^{p-1}(S^2L). 
 \]
 Clearly, the condition  $\gamma_L=0$ is equivalent to $\gamma_L^t=0$,

 \vskip .2cm
 
 \begin{prop}
 The system of maps $\gamma_L^t$, $L\in \SVect^\fd$,
  is a natural transformation of polynomial functors. 
 \end{prop}
 
 \noindent{\sl Proof:} To understand the nature of $\gamma_L$ and $\gamma_L^t$,
 we start from any representation $B$ of $\men$
 and  derive  from \eqref {eq:rels-in-g}
 the formulas  describing the action of $z\in L^*$ on
 $V_B^{\leq 2}$: 
 \be\label{eq:basic-form}
 \begin{gathered}
 z \ket{b} = 0,   \quad z y\ket{b}  = (-1)^{|y|\cdot |z|} yz \ket{b} + [z,y] \ket{b} = 
 -2 \ket{ e_{zy}(b)},
 \\
 z y_1y_2\ket{b} = (-1)^{|y|\cdot |z_1|} y_1 zy_2\ket{b} + [z,y_1] y_2 \ket{b} =
 \\
 = (-2) (-1)^{|y|\cdot |z_1|} y_1\ket{e_{z,y_2}(b)} -2
 e_{z, y_1} \bigl(y_2\ket{b}\bigr) =
 \\
 = (-2) (-1)^{|y|\cdot |z_1|}  y_1\ket{e_{z,y_2}(b)}  -2 e_{z,y_1}(y_2) \ket{b} + 
 (-2)
 (-1)^{ (|z|+|y_1|)\cdot |y_2|} y_2\ket {e_{z,y_1}(b)}
 \end{gathered}
 \ee
 for any  $y, y_1, y_2\in L$, $b\in B$. Here  the action of $e_{z,y_1}\in\men $ on products such
  as $y_2\ket{b}$ proceeds by the Leibniz rule. 
  
  \vskip .2cm
  
  We apply this to $B=\Lambda^{p-1}(S^2L)\otimes \CC_{-1/2}$. Then the action of
  an element of $\men$ like  $e_{zy}$ on $b\in \Lambda^{p-1}(S^2L)$
  should be the sum of the natural  action  on the tensor space
  and the additional term $-{1\over 2} \str(e_{zy})\Id$. Since
   $ \str(e_{zy}) = (z,y)$,
  we will have to add $-{1\over 2} (z,y) b$ to  
   $e_{zy}(b)$ and similarly in other cases. 
  
  \vskip .2cm
  
  Now look at $\gamma_L(z\otimes a) = zd\ket{a}$ with $a\in \Lambda^p(S^2L)$
 being the exterior product $a_1\wedge \cdots\wedge a_p$ with $a_i = \{y_i', y_i''\}
 \in S^2L$, so that $y'_i, y''_i\in L$. Applying  \eqref {eq:basic-form}  to $B$ as above
 as well as the formula for the Chevalley-Eilenberg differential $d$
 and the Moyal formula for the bracket in $\aen$ (Proposition \ref {prop:moyal}),
  we find that
 $\gamma_L(z\otimes a)$ is  a  linear combination of  ``monomial'' terms each of which
 is a product of:
 \begin{itemize}
 \item[(1)] The  scalar product of $z$ with some element of the set
 $S= \{y_1,\cdots, y'_p, y''_1, \cdots, y''_p\}$ (only one such scalar product can appear in
 each term).
 
 \item[(2)]  The product of two elements in S considered as an element of $U$.
 
 \item[(3)]  The wedge product of $p-1$ anticommutators of pairs of elements of $S$.

 \end{itemize}
  The coefficients in this linear combination will be independent of $L$. The fact that the scalar product in (1) appears in each term only once 
   implies that $\gamma_L^t$, written in the diagrammatic form
   (see, e.g., \cite [\S3.2]{SS-stable}), will not contain any  contractions $L\otimes L^*\to \CC$ or insertions
   $\CC\to L\otimes L^*$ and so  represents
    a morphism of polynomial functors. \qed
    
    \vskip .2cm
    
 Therefore (Remark \ref {rem:vanishing-even}) it suffices to verify vanishing of $\gamma_L$
 in the case when $L$ is a purely even space. We assume this from now on.
 
 \paragraph{$\gen$-linearity for $L$ even.}
  Let $L$ be purely even. Then $S^2(L)$ is spanned by perfect squares $y^2$,
  $y\in L$.  By Proposition \ref {prop:verma-lin} it therefore suffices to prove that
  \[
  z d\ket{ y_1^2\wedge\cdots \wedge y_p^2} = 0, \quad \forall \,\,p\geq 1, \,\,\,  y_1,\cdots, y_p\in L. 
  \]
 The formulas \eqref{eq:basic-form} imply for any $b\in \Lambda^{p-1}(S^2L)$:
 \be\label{eq:F1-F2}
 \begin{gathered}
 z[y_1, y_2] \ket{b}  =  2 (z,y_1) y_2\ket{b} -2  (z,y_2) y_1\ket{b},
 \\
 zy^2\ket{b} = -4y \ket{e_{zy}(b)} -2 (z,y) y\ket{b}. 
 \end{gathered}
 \ee
 Further, the Moyal formula of Proposition \ref{prop:moyal} (or an even easier direct
 calculation) gives the commutation law of perfect squares in $\aen$ in the form
 \[
 [y_1^2, y_2^2] = 2 [y_1, y_2] \cdot \{y_1, y_2\}. 
 \]
 This gives
 \be\label{eq:d(a)bef-z}
 \begin{gathered}
 d\ket{y_1^2\wedge \cdots \wedge y_p^2} \,\,  = \,\, \sum_{i=1}^p (-1)^{i-1} y_i^2
 \ket{y_1^2 \wedge\cdots \wedge \wh {y_i^2}\wedge\cdots \wedge y_p^2} -
 \\
 + 2 \sum_{i<j} (-1)^{i+j} [y_i, y_j] \,
 \ket{ \{y_i, y_j\} \wedge y_1^2 \wedge\cdots\wedge \wh {y_i^2} \wedge\cdots\wedge
 \wh {y_j^2}
 \wedge\cdots\wedge y_p^2}.
 \end{gathered}
  \ee
 Call the summands in the first sum in the RHS the {\em summands of the first kind}
 and those in the second sum the {\em summands of the second kind}. 
 By the first formula in \eqref {eq:F1-F2}, the action of $z\in L^*$ on a summand
 of the second kind does not affect the expression in the Dirac ``ket'', giving
 \be\label{eq:z-second}
 \begin{gathered}
 z[y_i, y_j] \ket{ \{y_i, y_j\} \wedge y_1^2 \wedge\cdots\wedge \wh{ y_i^2}
  \wedge\cdots\wedge \wh{y_j^2}
 \wedge\cdots\wedge y_p^2} = 
 \\
 =2 \bigl((z,y_i) y_j  - (z,y_j) y_i\bigr)\ket{ \{y_i, y_j\} \wedge y_1^2 \wedge\cdots\wedge\wh{ y_i^2}
  \wedge\cdots\wedge \wh{y_j^2}
 \wedge\cdots\wedge y_p^2} 
 \end{gathered} 
 \ee
 As for summands of the first kind,  the second formula in  \eqref {eq:F1-F2}
 gives
 \be\label{eq:F3}
 \begin{gathered}
 z y_i^2  \ket{y_1^2 \wedge\cdots \wedge \wh {y_i^2}\wedge\cdots \wedge y_p^2}\,=\,
 \\
 = \, -4 y_i \ket {e_{z,y_i}(y_1^2 \wedge\cdots\wedge \wh {y_i^2} \wedge\cdots\wedge
 y_p^2)}  -2  (z,y_i) y_i  
 \ket{y_1^2 \wedge\cdots \wedge \wh {y_i^2}\wedge\cdots \wedge y_p^2}.
 \end{gathered}
 \ee
 Here in the action of $e_{z,y_i}$ on the wedge product we use the Leibniz rule and add
 an extra term $-{1\over 2} \tr (e_{z,y_i})\cdot\Id  = - {1\over 2} (z,y_i)\cdot \Id$
 to account for tensoring with $\CC_{-1/2}$.
 This extra term will cancel the second term in the RHS of \eqref{eq:F3}. 
 
 \vskip .2cm
 
 Thus the action of $z$ on the sum of the summands of the first kind is
 \[
 4 \sum_{i=1}^p  (-1)^i y_i \ket {e_{z,y_i}(y_1^2 \wedge\cdots\wedge \wh {y_i^2} \wedge\cdots\wedge
 y_p^2)}
 \]
 where the action of $e_{z, y_i}$ is now without the twist by $\CC_{-1/2}$. 
 To expand this,  we note the
 general formula
 \[
 e_{zy}(y'{}^2) \,=\, \{ y', e_{zy}(y')\}  \, = \, (z,y') \{y,y'\}, 
 \]
and so by the Leibniz rule  we get 
 \[
 \begin{gathered}
 zy_i^2 \ket{y_1^2 \wedge\cdots \wedge \wh {y_i^2}\wedge\cdots \wedge y_p^2}\,\, = 
\\ 
= \,\, 4 \sum_{j\neq i} \eps_{ij} \cdot (z,y_j) y_i  \ket{ \{y_i, y_j\} \wedge y_1^2 \wedge\cdots\wedge \wh {y_i^2} \wedge\cdots\wedge
 \wh {y_j^2}
 \wedge\cdots\wedge y_p^2},
\end{gathered} 
 \]
 where $\eps_{ij} = (-1)^j $ for $i<j$ and $\eps_{ij}=(-1)^{j+1}$ for $i>j$. 
 Now \eqref {eq:z-second} shows that  the result of applying $z$ to the sum of summands of the first kind
 cancels with the result of applying it to the sum of summands of the second kind.
 That is, $zd\ket{y_1^2\wedge\cdots\wedge y_p^2}=0$.
So $d$ is $\gen$-linear, thus proving Theorem   \ref{thm:linearity}. 

 \vfill\eject

 \section{ Tensor categories} \label{sec:tensor}
 
 The constructions in this paper  are functorial in the vector superspace $L=\CC^{m|n}$. In fact, one can take for $L$ an
 object of a more general  category.
 
 \paragraph{ Bilinear pairs in a symmetric monoidal category and the Lie algebra $\o(2L+1)$. }\label{par:bilpairs}
  Let $(\Tc, \otimes)$ be a $\CC$-linear abelian symmetric monoidal category with 
  unit
 object $\b1$. We assume that $\Tc$ has countable direct sums and that $\otimes$ is exact in both arguments. 
 Our assumptions imply that we can form Schur functors $\Sigma^\alpha$  in $\Tc$, in particular, symmetric and exterior powers
 $S^p, \Lambda^p$.

 \vskip .2cm
 
 By a {\em bilinear pair} in $\Tc$ we will mean a datum of two objects $L, L_*\in\Tc$ and a pairing
 $\omega: L_* \otimes L\to\b1$. We make no assumptions as to (non)degeneracy of $\omega$; for example
 $\omega=0$ is allowed. In particular, 
   we do not require that $L_*$ is the categorical dual of $L$, although that may be the case. 

 \begin{prop}\label{prop:o(2L+1)}
 Given a bilinear pair $(L, L_*)$, the object
 \[
\o(2L+1) =  \gen(L, L_*) \,\,=\,\, \Lambda^2( L) \,\, \oplus \,\,  \  L\,\, \oplus  \,\, (  L\otimes   L_*)
 \,\, \oplus\,\,   L_* \,\, \oplus\,\, \Lambda^2 ( L_*)
 \]
 has a unique structure of a Lie algebra in $\Tc$ with the following properties: 
 
\begin{itemize}
\item[(1)]  The subobjects $\ol\nen = L\oplus \Lambda^2(L)$ and $\nen = L_* \otimes \Lambda^2(L_*)$  are Lie subalgebras with the only nonzero components of their brackets
being the wedge products $L\otimes L\to\Lambda^2(L)$ and $L_*\otimes L_* \to\Lambda^2(L_*)$. 

\item[(2)]  The subobject $\men = L\otimes L_*$ is a Lie subalgebra whose bracket comes from the associative algebra
structure defined by
\[
\men\otimes\men = L\otimes L_* \otimes L\otimes L_* \buildrel \Id\otimes\omega\otimes\Id \over \lra L\otimes L^* = \men
\]

\item[(3)]  The bracket of $\men = L\otimes L_*$ and $L\subset \ol\nen $ is the morphism $\Id\otimes\omega:
L\otimes L_*\otimes L \to L$

\item[(4)] The bracket of $\men$ and $L_*\subset \nen$ is the negative of the 
morphism 
\[
L\otimes L_* \otimes L_* \buildrel  R\otimes \Id \over
\lra L_*\otimes L\otimes L_* \buildrel \omega\otimes\Id\over\lra L_*,
\]
 where $R: L\otimes L_*\to L_*\otimes L$ is
the symmetry. 

\item[(5)]  The bracket of $L\subset\ol\nen$ and $L_*\subset \nen$ is twice  the identity map $L\otimes L_* \to L\otimes L_*=\men$. 
\end{itemize} 
 \end{prop}
The proposition can be either proved directly or deduced from the universal example in the next paragraph. 
  The reason for the  the notation $\o(2L+1)$   (slightly abusive as $L_*$ is not included) is the following
 
 \begin{ex}\label{ex:2C^{m|n}+1}
 Let $\Tc$ be the category of vector superspaces,  $L=\CC^{m|n}$, $L_*=L^*$ be the dual space and
 $\omega = \tr: L^*\otimes L\to \CC$  be the canonical pairing. Then Proposition \ref {prop:o(2L+1)} is true and
  $\o(2L+1)$ is the Lie superalgebra
 $\gen = \osp(2m+1|2n)$. The conditions (1-5) of the proposition are just the relations \eqref {eq:rels-in-g} in $\gen$. 
 \end{ex}
  
 \paragraph{The universal example: passing to infinity.} 
 Here is  a very important example of a bilinear pair
  \cite{DPS, SS-stable}. 
 Consider the Lie algebra $\mathfrak{sl}(\infty)$ defined as the direct limit of Lie algebras 
 $\mathfrak{sl}(m)$ with respect to the standard embeddings $\mathfrak{sl}(m)\to\mathfrak{sl}(m+1) $.
 Alternatively,  let $\mathbb L$ and $\mathbb L_*$ be  countably dimensional $\CC$-vector spaces equipped  with a
   pairing $\tr:\mathbb L_*\otimes \mathbb L\to\mathbb C$  which is nondegenerate, i.e., such that  its kernel in each of the arguments vanishes. It is easy to see \cite{mackey}  that a datum of $(\LL, \LL_*, \tr)$ is defined uniquely up to isomorphism, so we can
   assume that $\LL=\CC^\oo$ with the  basis $\{e_i\}_{i\geq 1}$ and $\LL_*$ is its restricted dual, spanned by the dual basis
   $\{e^i\}_{i\geq 1}$. 
 
 \vskip .2cm
 
  Then $\mathfrak{sl}(\infty)$ can be identified with the Lie algebra of all
 $X\in \End_{\mathbb C}(\mathbb L)$ such that $ X(\mathbb L)$ is finite dimensional, 
 $X^t$ (given by the transposed matrix in our basis) takes $\LL_*$ to a finite-dimensional subspace in $\LL_*$
 and 
  $\tr (X)=0$. 
  
   \vskip .2cm
 By construction,  $\mathbb L$ and $\mathbb L_*$ are  $\mathfrak{sl}(\infty)$-modules, and $\ssl(\oo)$ can be identified naturally with the subalgebra of $\gl(\infty)\simeq \mathbb L\otimes \mathbb L_*$ given by the kernel of the pairing $\mathbb L\otimes \mathbb L_*\to\mathbb C$.
  Let $\mathbb T_\fsl$ denote the abelian category of $\fsl(\infty)$-modules generated by the tensor powers of $\mathbb L$ and $\mathbb L_*$. This is a symmetric monoidal category of the kind considered in  \S \ref {par:bilpairs}.
   It  was studied in \cite{DPS} and in \cite{SS-stable}. In particular, the following universality result was
  proved in \cite{SS-stable}.
  
  \begin{thm}\label{thm:SS}
  Let $(L, L_*, \omega)$ be a bilinear pair in a  symmetric  monoidal category $\Tc$ as \S \ref {par:bilpairs}. Then there
  is a unique, up to isomorphism,  left-exact symmetric monoidal functor $F:\mathbb T_\fsl\to \mathcal T$ such that 
  $F(\mathbb L)=L$, $F(\mathbb L_*)=L_*$ and $F(\tr)=\omega$.  
 \end{thm}
 
 \noindent{\sl Proof:}  More precisely, \cite[Th.3.4.2]{SS-stable} proves this kind of statement for
 the category $\Rep(GL)$ of representations of the group $GL=\varinjlim\,  GL(m)$ generated, as an abelian
 symmetric monoidal category, by the standard representation $\LL=\CC^\oo$ and its restricted dual $\LL_*$.  
 But it is immediate that $\Rep(GL)$ is equivalent to $\TT_\ssl$, as the morphisms between the  tensor powers
 $\LL^{\otimes p} \otimes \LL_*^{\otimes q}$  are the same in both categories. 
 \qed 
 
 \begin{ex}
 The object $\o(2\LL+1) = \gen(\LL, \LL_*)$ in $\TT_\ssl$ is nothing but  the Lie algebra $\o(2\oo+1) = \varinjlim \o(2m+1)$,
 the inductive limit with respect to the standard embeddings and so Proposition \ref {prop:o(2L+1)} holds in this
 case. 
 This follows from Example \ref  {ex:2C^{m|n}+1} which established the identification for each finite $m$. 
 Applying Theorem \ref{thm:SS}, we get a calculation-free proof of  Proposition \ref {prop:o(2L+1)} for a general bilinear pair. 
    \end{ex}
    
    \paragraph{The spinor module over $\o(2L+1)$ and its BGG-type  resolution.}   Now,  categorical generalizations
    of all the other constructions and results of this paper can be obtained in a straightforward way. 
    
    \vskip .2cm
    
Let  $(L, L_*, \omega)$ be a bilinear pair
    in $(\Tc, \otimes)$ as before and let $\gen=\o(2L+1)$.
    For any $\eps\in \CC$ we denote by $\b1_\eps$ the representation of the Lie algebra
    $\men = L\otimes L_*$ on the object $\b1$ given by the character $\eps\cdot \omega: \men\to\CC$. 
    
   \begin{prop}
   The object $\Lambda^\bullet(L) = \bigoplus_p \Lambda^p(L)$ has a structure of $\o(2L+1)$-module defined uniquely by the properties:
   \begin {itemize}
   \item[(1)] The action of $L\subset \ol\nen\subset\gen $ is given by the wedge product 
   (antisymmetrization) $L\otimes\Lambda^p(L)\to \Lambda^{p+1}(L)$.
   
   \item[(2)] The action of $L_*\subset \nen\subset \gen$ is given by the contraction 
   $L_* \otimes\Lambda^p(L)\to \Lambda^{p-1}(L)$. 
   \end{itemize}
   \end{prop}
   
   \noindent {\sl Proof:} One can either prove this directly, or use Proposition \ref {prop:spin/weil} to deduce it for
   the universal bilinear pair $(\LL, \LL_*, \tr)$ in $\TT_\ssl$ and then  apply Theorem \ref {thm:SS}. \qed
    
    \vskip .2cm
    
    Let us call $\Lambda^\bullet(L)$ the {\em spinor module} over $\gen = \o(2L+1)$.     A BGG-type resolution 
    $\Vc_\Tc^\bullet = \Cc_\Tc^\bullet\otimes \b1_{-1/2}$ of
     $\Lambda^\bullet(L)$ can be obtained by repeating the constructions of  \S \ref {const-resul}\ref {par:exp-con}
     in our categorical context. We leave the details to the reader, indicating just one point.
     
     \vskip .2cm
     
     After $\Vc^\bullet_\Tc$ has been constructed, the fact that  it is indeed a resolution of $\Lambda^\bullet(L)$
     follows from Theorem \ref{thm:SS}. More precisely,   it suffices to verify this fact  for $\Cc^\bullet_\Tc$, as tensoring
     with $\b1_{-1/2}$ does not affect it. That is, we need to prove that the augmented complex $\wt\Cc^\bullet_\Tc = \on{Cone}
     \{ \Cc_\Tc^\bullet\to\Lambda^\bullet(L)\}$ is exact. 
     For this,  we first consider the universal case $\Tc=\TT_\ssl$, $L=\LL$ etc.
     The exactness 
   of   $\wt \Cc^\bullet_{\TT_\ssl} $  
      follows from Theorems \ref {thm:exactness}  and \ref{thm:linearity} by passing to the inductive limit.

          \vskip .2cm
      
      In the general case we use   the   functor $F: \TT_\ssl\to\Tc$ from Theorem \ref{thm:SS} taking $\LL\to L$ etc. 
      Every term of $\wt\Cc_{\TT_\ssl}^\bullet$ is a direct sum of some $\Sigma^{\alpha}\mathbb L$ and $F(\Sigma^{\alpha}\mathbb L)=\Sigma^\alpha(L)$.
      Now, in $\mathbb T_{\fsl}$ we have the orthogonality relations:  
      \[
      \Hom(\Sigma^\alpha(\mathbb L),\Sigma^\beta(\mathbb L)) = \begin{cases} \CC, & \text{ if } \alpha=\beta,
      \\
      0, & \text {otherwise}.
      \end{cases}
      \]
    Since the differentials in $\wt\Cc_{\TT_\ssl}^\bullet$
     are sums of morphisms of Schur functors, these relations imply that $\wt\Cc^\bullet_\TT $ splits into direct sum of complexes
      $\wt\Cc^\bullet_\TT(\alpha) \otimes \Sigma^\alpha(\LL)$, where
       $\wt\Cc^\bullet_\TT(\alpha) = \Hom(\Sigma^\alpha(\LL), \wt\Cc^\bullet_\TT)$  is the complex of vector spaces formed
     by the  multiplicities of $\Sigma^\alpha(\LL)$ in $\wt\Cc^\bullet_\TT$. Therefore 
     exactness of $\wt\Cc^\bullet_\TT$ implies the exactness of each 
  $\wt\Cc^\bullet_\TT(\alpha )$ as a complex of vector spaces.   Because $F$ is a $\CC$-linear functor, we conclude that $\wt\Cc^\bullet_\Tc = F(\wt\Cc^\bullet_\TT) \=
  \bigoplus_\alpha \wt\Cc^\bullet_\TT(\alpha) \otimes \Sigma^\alpha (L)$ has each summand exact and so is exact itself.

 \vfill\eject

	\vskip 1cm

M.K.: Kavli IPMU (WPI), UTIAS, The University of Tokyo,
5-1-5 Kashiwanoha, Kashiwa, Chiba, 277-8583 Japan. 
Email: {\tt mikhail.kapranov@protonmail.com} 

\vskip .2cm

V.S.: Department of Mathemtics, University of California, Berkeley CA
94720. Email: {\tt serganov@math.berkeley.edu}

\end{document}